\RequirePackage{fix-cm}
\documentclass[smallextended,nospthms]{svjour3}
\usepackage[T1]{fontenc}
\usepackage{lmodern}
\usepackage{enumerate}
\usepackage{mathrsfs,graphicx,bm,dsfont}
\usepackage{xcolor}
\usepackage{amsthm,amsmath,amsfonts,amssymb}
\usepackage{algorithm,algorithmic}
\usepackage{mathtools,needspace,cite,etoolbox,xurl}
\usepackage[nopatch=footnote]{microtype}
\usepackage[colorlinks=true,linkcolor=blue,citecolor=blue,urlcolor=blue,
  anchorcolor=blue,breaklinks=true]{hyperref}
\providecommand{\doi}[1]{\url{https://doi.org/#1}}
\DeclareMathSymbol{\Omega}{\mathalpha}{operators}{"0A}
\DeclareMathSymbol{\Theta}{\mathalpha}{operators}{"02}
\allowdisplaybreaks[3]
\journalname{Mathematical Programming}

\theoremstyle{plain}
\newtheorem{theorem}{Theorem}[section]
\newtheorem{lemma}[theorem]{Lemma}
\newtheorem{proposition}[theorem]{Proposition}
\newtheorem{corollary}[theorem]{Corollary}

\theoremstyle{definition}
\newtheorem{definition}[theorem]{Definition}
\newtheorem{assumption}{Assumption}[section]

\theoremstyle{remark}
\newtheorem{remark}[theorem]{Remark}

\newcommand{\RR}{\mathbb{R}}

\newcommand{\dist}{\mathrm{dist}}
\newcommand{\diam}{\mathrm{diam}}
\newcommand{\proj}{\mathrm{Proj}}
\newcommand{\cX}{\mathcal{X}}

\renewcommand{\hat}{\widehat}

\newcommand{\bx}{\mathbf{x}}
\newcommand{\bc}{\mathbf{c}}

\newcommand{\by}{\mathbf{y}}
\newcommand{\bw}{\mathbf{w}}
\newcommand{\bn}{\mathbf{n}}
\newcommand{\bg}{\mathbf{g}}
\newcommand{\bz}{\mathbf{z}}
\newcommand{\bv}{\mathbf{v}}
\newcommand{\bu}{\mathbf{u}}
\newcommand{\bp}{\mathbf{p}}

\newcommand{\bq}{\mathbf{q}}
\newcommand{\ba}{\mathbf{a}}
\newcommand{\br}{\mathbf{r}}

\newcommand{\ball}{\mathbb{B}}

\newcommand{\unit}{\mathbb{S}^{d-1}}

\DeclareMathOperator{\Lip}{Lip}

\numberwithin{equation}{section}

\hypersetup{
  pdftitle={Deterministic Oracle Complexity for Goldstein Stationarity in Nonsmooth Nonconvex Optimization with Constraints},
  pdfauthor={Yixi Ding and Ya-xiang Yuan},
  pdfsubject={Deterministic oracle complexity lower bounds for constrained nonsmooth nonconvex optimization},
  pdfkeywords={Goldstein stationarity, deterministic oracle complexity, nonsmooth nonconvex optimization, set constraints, functional constraints},
  bookmarksnumbered=true,
  hypertexnames=false
}

\title{Deterministic Oracle Complexity for Goldstein Stationarity in Nonsmooth Nonconvex Optimization with Constraints}
\titlerunning{Deterministic Oracle Complexity for Constrained Goldstein Stationarity}
\author{Yixi Ding \and Ya-xiang Yuan}
\authorrunning{Yixi Ding and Ya-xiang Yuan}
\institute{Yixi Ding (corresponding author) \at
State Key Laboratory of Mathematical Sciences, ICMSEC,
Academy of Mathematics and Systems Science,
Chinese Academy of Sciences, Beijing 100190, China;\\
University of Chinese Academy of Sciences, Beijing 100049, China\\
\email{dingyixi@amss.ac.cn}\\
ORCID: \url{https://orcid.org/0009-0009-9290-1909}
\and
Ya-xiang Yuan \at
State Key Laboratory of Mathematical Sciences, ICMSEC,
Academy of Mathematics and Systems Science,
Chinese Academy of Sciences, Beijing 100190, China\\
\email{yyx@lsec.cc.ac.cn}}
\date{}

\begin{document}
\maketitle
\begin{abstract}
We consider the deterministic oracle complexity of finding a $(\delta,\epsilon)$-Goldstein stationary point of a Lipschitz function under constraints. Our main result concerns products of convex bodies. Each block contains an inner ball and lies in a concentric outer ball. The inner radii have a common positive lower bound, and each outer radius is at most a fixed multiple of the corresponding inner radius. This class includes common constraints such as balls and boxes. On each prescribed feasible set of dimension $d$ in this class, every deterministic algorithm using a local oracle needs $\exp(\Omega(d))$ objective-oracle calls in the worst case at fixed, sufficiently small accuracies $\delta$ and $\epsilon$. The objective Lipschitz and initial objective-value gap bounds are fixed, and each constructed function has a global minimizer in the interior of the feasible set. For broader classes of constraints, supplementary results show that every deterministic algorithm needs $\Omega(\epsilon^{-2})$ objective-oracle calls in the worst case when $d\gtrsim\epsilon^{-2}$ and $\delta\lesssim\epsilon$. The constants depend on the interior-ball radius and, for inequalities, the local Lipschitz and Slater bounds. These supplementary results cover closed sets that are star-shaped with respect to a ball and convex coercive inequalities with a Slater point, and hold even for convex objectives.
\end{abstract}

\section{Introduction}
\label{sec:intro}

Constrained nonsmooth nonconvex problems arise in many application areas, including machine learning, operations research, and statistical estimation. Their analysis is difficult because the objective is not required to be smooth or convex, and first-order necessary conditions involve both the objective subdifferential and the normal cone to the feasible set. First-order methods are widely used for such problems \cite{SFOM, OCSS, FOM}. In this paper, we consider the problem
\begin{equation}
\label{eq:problem}
\min_{\bx \in \cX} f(\bx),
\end{equation}
where $d$ is the dimension, $f:\mathbb{R}^d\to\mathbb{R}$ is Lipschitz continuous, and $\cX\subseteq\mathbb{R}^d$ is closed. For unconstrained problems ($\cX=\mathbb{R}^d$) with Lipschitz continuous gradients and a bounded initial objective-value gap $f(\bx_0)-\inf_{\bx\in\RR^d}f(\bx)$, measured at the initial point $\bx_0$, deterministic gradient algorithms can find a point $\bx$ satisfying $\|\nabla f(\bx)\|\le\epsilon$ with dimension-free $\mathcal{O}(\epsilon^{-2})$ gradient-oracle complexity \cite{Nesterov18, Nesterov12}.  In the nonsmooth Lipschitz setting, however, even computing a point $\bx$ satisfying $\dist(\mathbf 0,\partial f(\bx))\le\epsilon$, where $\partial f(\bx)$ denotes the Clarke subdifferential of $f$ at $\bx$, is impossible for any finite-time randomized or deterministic algorithm interacting with a local oracle in full generality \cite{nogas, KS33, TS21}.  This obstruction motivates approximate stationarity notions that aggregate first-order information over a neighborhood of radius $\delta$. \textit{$(\delta,\epsilon)$-Goldstein stationarity} is one such notion: instead of using only $\partial f(\bx)$ at one point $\bx$, it tests whether the convex hull of the union of Clarke subdifferentials $\partial f(\by)$ over points $\by$ within distance $\delta$ of $\bx$ contains a vector of norm at most $\epsilon$ \cite{Goldstein, GSL, Davis22}.  In unconstrained Lipschitz optimization, randomized algorithms can find $(\delta,\epsilon)$-Goldstein stationary points with dimension-free oracle complexity bounds \cite{CFS, Davis22, TS22}.  By contrast, deterministic algorithms cannot in general find such points within the same dimension-free finite number of oracle calls \cite{nogas, DNN, KS22}.

In the constrained problem, the \textit{first-order optimality condition} also depends on the feasible-set geometry. At a feasible point $\bx$, this condition takes the form \cite[Theorem 8.15]{V44}:
\begin{equation}
\label{eq:first_order_optimality_intro}
    \mathbf 0\in \partial f(\bx)+N_\cX(\bx).
\end{equation}
This leads to the question of whether deterministic algorithms can compute approximate stationary points and, if they can, with what oracle complexity. The algorithms considered here interact with a \textit{local oracle} in the sense of \cite{nogas}: at a queried point $\bx$, the oracle returns a function agreeing with $f$ on some neighborhood of $\bx$, without specifying the neighborhood radius. Such a response provides the function value, the entire Clarke subdifferential, and higher-order derivatives when they exist. The algorithm may be given the entire feasible set or explicit constraint function $g$. When considering the set-constrained problem \eqref{eq:problem}, although this oracle does not return normal vectors, the normal cone is determined by the feasible set and enters the optimality condition \eqref{eq:first_order_optimality_intro}. A feasible point $\bx$ is \textit{$(\delta,\epsilon)$-Goldstein stationary} if the following \textit{stationarity measure} is at most $\epsilon$:
\begin{equation}
\label{eq:constrained_goldstein_measure_intro}
    \dist\biggl(\mathbf 0,\operatorname{conv}\Bigl( \bigcup_{\by \in \ball_\delta(\bx) \cap \cX}
    \bigl(\partial f(\by) + N_\cX(\by)\bigr) \Bigr)\biggr).
\end{equation}
To prove a lower bound on the number of objective-oracle calls under set constraints, the construction must produce an objective consistent with all local responses for which the stationarity measure at the algorithm's returned feasible point exceeds $\epsilon$. Ensuring this inequality requires more than controlling the objective subdifferential terms: the construction must also ensure that adding normal vectors throughout the $\delta$-neighborhood does not bring the convexified set in \eqref{eq:constrained_goldstein_measure_intro} close to the origin.

An example illustrates this point. Suppose that a hard construction from the unconstrained setting has a subgradient $\mathbf e_1\in\partial f(\mathbf 0)$.  If the same objective function $f$ is placed over the nonnegative orthant $\cX=\mathbb R_+^d$ and the algorithm queries $\mathbf 0$, then
\[
    N_\cX(\mathbf 0)=\mathbb R_-^d,
    \qquad -\mathbf e_1\in N_\cX(\mathbf 0).
\]
Consequently,
\[
    \mathbf 0\in \partial f(\mathbf 0)+N_\cX(\mathbf 0).
\]
This illustrates why hard functions from the unconstrained setting cannot be simply transplanted to the constrained one. A central technical issue for the set-constrained problem is therefore to construct the objective function $f$ while accounting for the normal cones of the feasible set.

The main result concerns products of convex bodies. Each block contains a ball of radius at least $r_{\mathrm{in}}>0$ and is contained in a concentric ball whose radius is at most $K$ times the inner radius. This class includes Euclidean balls, boxes, and products of balls and ellipsoids with uniformly bounded axis ratios. For each prescribed feasible set in this class, Theorem~\ref{thm:exp_set} gives an $\exp(\Omega(d))$ worst-case lower bound on the number of objective-oracle calls needed by any deterministic algorithm to find a point whose stationarity measure in \eqref{eq:constrained_goldstein_measure_intro} is at most $\epsilon$. Here $L_f$ bounds the objective Lipschitz constant, and $\Delta_f$ bounds the initial objective-value gap $f(\bc)-\inf_{\bx\in\cX}f(\bx)$ at the block center $\bc$. The bounds $L_f,\Delta_f$, the common interior radius $r_{\mathrm{in}}$, and the ratio bound $K$ are independent of $d$. The conclusion holds at fixed, sufficiently small stationarity accuracy $\epsilon$ and neighborhood radius $\delta$, as specified in Corollary~\ref{cor:exp_fixed_accuracy}.

The exponential construction starts from a query lower bound for approximate Brouwer fixed points \cite{GKP22}. A Lipschitz vector field derived from the Brouwer map prescribes the gradients of the affine oracle responses. These gradients vary across query points. An adaptation of the finite-set interpolation construction of Ding and Yuan \cite[Lemma 3.2]{DY26} realizes all these responses as a single Lipschitz objective. In the sampling neighborhood of the returned point, all objective subgradients have a uniformly positive projection onto a common direction, and all normal vectors at feasible points have a nonnegative projection onto that direction. The constructed function has a global minimizer in the interior of the feasible set. The argument also gives an exponential lower bound on the number of objective-oracle calls needed to find a point satisfying the generalized Goldstein criterion of Liu et al.\ \cite{zero}, with a prescribed projection parameter $\gamma$ satisfying the condition in Corollary~\ref{cor:ggsp_exp}.

For more general constraints, supplementary bounds in terms of the stationarity accuracy follow from a resisting-oracle construction (see, e.g., \cite{nemirovski1983problem, Nesterov18} for foundational theory of resisting oracles, and \cite{KS22, DNN, nogas} for recent applications in unconstrained nonsmooth optimization). The objective is constructed as the pointwise maximum of affine functions with mutually orthogonal gradients and a radial term centered at an interior reference point. The affine components preserve earlier local responses, while the radial term controls the effect of normal vectors near the boundary. The argument applies to closed feasible sets that are star-shaped with respect to an interior ball, including all convex bodies with nonempty interior. For fixed interior-ball radius, objective Lipschitz bound, and function-value gap bound, every deterministic algorithm needs $\Omega(\epsilon^{-2})$ objective-oracle calls in the worst case when $d\gtrsim\epsilon^{-2}$ and $\delta\lesssim\epsilon$. The same construction gives a Goldstein KKT (GKKT) lower bound of this order for every prescribed finite-valued convex coercive constraint function with a Slater point. In both supplementary results, the constructed function is convex and has a global minimizer in the interior of the feasible set. The $\epsilon^{-2}$ order is already obtainable on balls from classical convex minimization lower bounds; the supplementary results extend the construction to each prescribed constraint in the stated broader classes.

The constrained criteria used in the literature differ in where objective subgradients and normal vectors are sampled and in how projection steps enter the definition. Section~\ref{sec:preliminaries} compares \eqref{eq:constrained_goldstein_measure_intro} with the criteria of Liu et al.\ \cite{zero} and Buzaglo and Kornowski \cite{BK26}. Counterexamples show that the joint neighborhood criterion is generally incomparable with either criterion; additional geometric or smoothness conditions give quantitative implications in both directions.

\subsection{Related work}
\label{subsec:related_work}
With the objective Lipschitz bound, initial objective-value gap bound, and success probability fixed, a randomized value--subgradient algorithm finds a $(\delta,\epsilon)$-Goldstein stationary point using $\widetilde{\mathcal O}(\delta^{-1}\epsilon^{-3})$ oracle calls, where the tilde suppresses logarithmic factors \cite[Theorem 6]{Davis22}. Related computational procedures approximate the relevant subdifferential information using gradient sampling or zeroth-order smoothing \cite{GSA, KS24, ZONN, zero, ISC}.

A separate line of work establishes lower bounds on the oracle complexity of deterministic algorithms for unconstrained nonsmooth problems. At fixed admissible accuracies, Jordan et al.\ \cite[Theorem 1]{DNN} prove an $\Omega(d)$ first-order oracle lower bound, and Tian and So \cite[Theorem 1]{nogas} obtain an $\Omega(d)$ bound for deterministic algorithms using a local oracle. These results separate deterministic complexity from the dimension-free randomized bounds. Kornowski and Shamir \cite[Theorem 3]{KS22} also establish an $\Omega(\epsilon^{-2})$ first-order lower bound for unconstrained convex objectives, even for randomized algorithms. Their proof combines a max-type function in orthogonal coordinates with a radial term. The precision-dependent construction in Section~\ref{sec:framework1} adapts this argument to control normal vectors throughout the Goldstein neighborhood and gives a joint lower bound in $d$, $\epsilon$, and $\delta$ for all deterministic algorithms using a local oracle. Its scope includes every prescribed closed feasible set that is star-shaped with respect to an interior ball, and the constructed functions remain convex. Classical local-oracle lower bounds for convex minimization \cite[Section 4.4.3, Proposition 1]{nemirovski1983problem} also imply this $\epsilon^{-2}$ order on a ball through the estimates in Proposition~\ref{prop:objective_gap} relating approximate stationarity to the objective-value gap. Kornowski et al.\ \cite[Theorem 2.1]{HMLG} establish an exponential-in-dimension oracle lower bound, even for randomized local algorithms. Their result concerns finding an iterate whose function-value gap relative to the minimum over a ball centered at that iterate is at most a fixed small fraction of the ball's radius.

With the objective Lipschitz and initial objective-value gap bounds fixed, Kornowski \cite[Theorem 3.1]{KFO26} establishes an $\Omega(\delta^{-1}\epsilon^{-3})$ first-order lower bound for unconstrained Goldstein stationarity. Its parameter range requires $\epsilon\lesssim\delta\lesssim\epsilon^{-1}$, sufficiently small $\epsilon$, and sufficiently large dimension. The formal statement concerns randomized zero-respecting algorithms, even when the oracle provides the entire Clarke subdifferential. This result complements the dimension-free randomized upper bound above.

More recently, Kornowski \cite[Theorem 3.1]{K26} and Ding and Yuan \cite[Theorem 2.2]{DY26} proved that, for unconstrained Lipschitz optimization, any deterministic algorithm using a local oracle requires $(1/\epsilon)^{\Omega(d)}$ oracle calls in the worst case to find a $(\delta,\epsilon)$-Goldstein stationary point. These bounds apply with the objective Lipschitz and initial objective-value gap bounds fixed, for sufficiently small positive $\delta,\epsilon$ and sufficiently large $d$. Ding and Yuan \cite[Theorem 2.3]{DY26} also give a deterministic value-only upper bound of the same exponential order in $d$ at fixed accuracies. The set-constrained exponential result additionally controls normal vectors at all feasible samples in the neighborhood.

For constrained nonsmooth nonconvex optimization, recent algorithmic work studies inequality and set constraints using first-order models, sampling schemes, or Goldstein-type stationarity criteria \cite{SFOM, OCSS, FOM, GSL}. Liu et al.\ \cite[Definitions 4.1--4.2]{zero} introduce generalized Goldstein stationarity through a projected gradient mapping with a prescribed step and give stochastic zeroth-order algorithms. Buzaglo and Kornowski \cite[Section 2]{BK26} use the ambient Goldstein subdifferential together with the normal cone at the current point. Their residual bounds Liu et al.'s fixed-step residual from above. With the objective Lipschitz, initial objective-value gap, and stochastic-gradient second-moment bounds fixed, their algorithm uses $\widetilde{\mathcal O}(\delta^{-1}\epsilon^{-3})$ calls to obtain an expected residual of at most $\epsilon$ \cite[Theorem 4.1]{BK26}. For the same residual on a ball, they give an $\Omega(\delta^{-1}\epsilon^{-3})$ stochastic-gradient lower bound under the joint parameter conditions of \cite[Corollary 8.23]{BK26}. This lower bound applies to randomized zero-respecting algorithms whose outputs are supported on coordinates appearing in the observed stochastic gradients. The measure in \eqref{eq:constrained_goldstein_measure_intro} uses normal vectors at nearby feasible samples. Propositions~\ref{prop:ggsp_comparison} and~\ref{prop:residual_comparison} compare these criteria.

Recent work treats further bounded-domain and composite models. Zhang et al.\ \cite{ZWZ26} prove a stochastic zeroth-order lower bound with queries and outputs restricted to a Euclidean ball; their target uses the ambient Goldstein subdifferential and contains no feasible-set normal cone. Li et al.\ \cite{LTY26} give randomized complexity upper bounds for proximal Goldstein stationarity in composite optimization, including convex set constraints. When the regularizer is the indicator function of $\cX$, their criterion coincides with Liu et al.'s fixed-step criterion \eqref{eq:ggsp_residual}.

\subsection{Organization}

The remainder of the paper is organized as follows. Section~\ref{sec:preliminaries} collects notation, stationarity concepts, comparisons between constrained criteria, and the oracle model. Section~\ref{sec:framework1} gives the precision-dependent lower bound under set constraints. Section~\ref{sec:exponential} establishes the main exponential lower bound and its consequence for the fixed-step criterion. Section~\ref{sec:inequality_constraints} gives the precision-dependent result under convex inequalities. Section~\ref{sec:discussion} concludes with open questions. The appendices contain the precision-dependent constructions, auxiliary estimates for the exponential results, and proofs of the criterion comparisons.

\section{Preliminaries and geometric background}
\label{sec:preliminaries}

\subsection{Notation and stationarity measures}

\paragraph{Notation.} 

The notation used in this paper is standard in nonsmooth optimization. Scalars are denoted by lowercase letters (e.g., $\alpha, \delta$), and vectors by boldface lowercase letters (e.g., $\bx, \bv$). Feasible constraint domains are denoted by calligraphic uppercase letters (e.g., $\mathcal{X}, \mathcal{S}$), and general sets (e.g., open neighborhoods) by standard uppercase letters (e.g., $U, V$). The ambient space is the Euclidean space $\mathbb{R}^d$ equipped with the standard inner product $\langle \cdot, \cdot \rangle$ and Euclidean norm $\|\cdot\|$. The standard basis vectors are $\mathbf{e}_1, \dots, \mathbf{e}_d$.

For a set $S \subseteq \mathbb{R}^d$, its interior, closure, and boundary are denoted by $\operatorname{int}(S)$, $\overline{S}$, and $\partial S$, respectively. The distance from a point $\bx$ to a set $S$ is $\operatorname{dist}(\bx, S) := \inf_{\by \in S} \|\bx - \by\|$. The convex hull of $S$ is denoted by $\operatorname{conv}(S)$. The closed Euclidean ball centered at $\bx$ with radius $r$ is denoted by $\ball_{r}(\bx)$, while the corresponding open ball is denoted by $B(\bx,r)$. For a Lipschitz continuous function $f$, its Lipschitz constant is denoted by $\operatorname{Lip}(f)$. For a finite set $S$, its cardinality is denoted by $\#S$. For $t\in\RR$, the symbols $\lfloor t\rfloor$ and $\lceil t\rceil$ denote the greatest integer not exceeding $t$ and the least integer not smaller than $t$, respectively. The symbols $\mathcal{O}$ and $\Omega$ denote asymptotic upper and lower bounds, respectively; the constants are independent of the parameters varied in the stated bound. 

\paragraph{Stationarity measures.}
The analysis focuses on locally Lipschitz continuous functions on $\mathbb{R}^d$. By Rademacher's theorem \cite{Clarke1990, V44}, such functions are differentiable almost everywhere. This allows us to define the Clarke subdifferential and its elements, called Clarke subgradients.

\begin{definition}[Clarke Subdifferential \cite{Clarke1990}]\label{def:Clarke_Subdifferential}
Let $f: \mathbb{R}^d \to \mathbb{R}$ be locally Lipschitz. The Clarke subdifferential of $f$ at $\bx$, denoted $\partial f(\bx)$, is the convex hull of all limit points of gradients:
\[
  \partial f(\bx) := \operatorname{conv} \left\{ \bg \mid \exists \bx_k \to \bx \text{ with } f \text{ differentiable at } \bx_k \text{ and } \nabla f(\bx_k) \to \bg \right\}.
\]
Any individual vector $\bg \in \partial f(\bx)$ is called a \textit{Clarke subgradient}.
\end{definition}

While the Clarke subdifferential captures local geometry, it is not by itself the most robust object for black-box finite-time analysis. We use a stationarity criterion that aggregates subdifferential information over a neighborhood $\ball_\delta(\bx)$.

\begin{definition}[Goldstein Subdifferential \cite{Goldstein}]
Given a sampling radius $\delta > 0$, the standard Goldstein $\delta$-subdifferential is the convex hull of the union of Clarke subdifferentials over the $\delta$-ball:
\[
  \partial_\delta f(\bx) := \operatorname{conv}\bigg( \bigcup_{\by \in \ball_\delta(\bx)} \partial f(\by) \bigg).
\]
\end{definition}

To define an appropriate measure of stationarity for constrained problems, the local geometry of the feasible set also needs to be taken into account. The standard normal cones from variational analysis are recalled below; see, for example, \cite{V44}.
\begin{definition}[Limiting Normal Cone]
\label{def:rw-normals}
Let $\cX \subseteq \mathbb{R}^d$ be a closed set and $\bar{\bx} \in \cX$. A vector $\bv$ is a \textit{regular normal} to $\cX$ at $\bar{\bx}$, written $\bv \in \hat{N}_\cX(\bar{\bx})$, if $\langle \bv, \bx - \bar{\bx} \rangle \leq o(\|\bx - \bar{\bx}\|)$ for all $\bx \in \cX$. A vector $\bv$ is a \textit{limiting normal} to $\cX$ at $\bar{\bx}$, written $\bv \in N_\cX(\bar{\bx})$, if there exist sequences $\bx^\nu \to \bar{\bx}$ in $\cX$ and $\bv^\nu \to \bv$ with $\bv^\nu \in \hat{N}_\cX(\bx^\nu)$.
\end{definition}

These normal cones connect classical variational analysis with the neigh\-bor\-hood-based stationarity notion used below. For a local optimum $\bx^* \in \cX$, a standard first-order necessary condition is
\begin{equation}\label{first-order necessary condition}
  \mathbf{0} \in \partial f(\bx^*) + N_\cX(\bx^*),
\end{equation}
as given, for example, by \cite[Theorem 8.15]{V44}. This exact condition \eqref{first-order necessary condition} motivates its neighborhood-based analogue and the corresponding approximate stationarity criterion.

\begin{definition}[Constrained Goldstein Stationarity]\label{def:constrained_goldstein}
\label{def:goldstein-stationarity}
Given a feasible point $\bx \in \cX$ and a sampling radius $\delta > 0$, the \textit{constrained Goldstein subdifferential} is defined by
\[
  \partial_\delta^\cX f(\bx) := \operatorname{conv}\bigg( \bigcup_{\by \in \ball_\delta(\bx) \cap \cX} \big(\partial f(\by) + N_\cX(\by)\big) \bigg).
\]
A point $\bx \in \cX$ is called a \textit{$(\delta,\epsilon)$-Goldstein stationary point} if the distance from the origin to its constrained Goldstein subdifferential is bounded by $\epsilon$:
\[
  \operatorname{dist}\big(\mathbf{0}, \; \partial_\delta^\cX f(\bx) \big) \le \epsilon.
\]
\end{definition}

Write $S_\delta(\bx)=\dist(\mathbf0,\partial_\delta^\cX f(\bx))$ for the set-constrained stationarity measure.

\paragraph{Generalized Goldstein stationarity.}
Liu et al.\ \cite[Definitions 4.1--4.2]{zero} introduce a projection criterion on compact convex feasible sets. The same formula defines a residual on any nonempty closed convex set $\cX$. For $\bx\in\cX$ and a fixed step $\gamma>0$, we denote the minimum in their criterion by
\begin{equation}\label{eq:ggsp_residual}
M_{\gamma,\delta}(\bx)
=\min_{\bv\in\partial_\delta f(\bx)}
\frac{\|\bx-\proj_{\cX}(\bx-\gamma \bv)\|}{\gamma}.
\end{equation}
Their $(\gamma,\delta,\epsilon)$-generalized Goldstein stationary point (GGSP) is a feasible point with $M_{\gamma,\delta}(\bx)\le\epsilon$. In \eqref{eq:ggsp_residual} the subgradient samples may come from the entire ambient ball. The normal vector associated with the projection is taken at the projected point. In Definition~\ref{def:constrained_goldstein} both the objective subgradient and the normal vector are evaluated at each feasible sampling point before convexification. 

For $\bv\in\partial_\delta f(\bx)$, write
\[
\bp_{\bv}=\proj_{\cX}(\bx-\gamma \bv),\qquad
\bw_{\bv}=\frac{\bx-\bp_{\bv}}{\gamma},\qquad
\bn_{\bv}=\frac{\bx-\gamma \bv-\bp_{\bv}}{\gamma}.
\]
The projection optimality condition is
\begin{equation}\label{eq:ggsp_projection}
\bw_{\bv}=\bv+\bn_{\bv},\qquad \bn_{\bv}\in N_{\cX}(\bp_{\bv}).
\end{equation}

The ambient Goldstein subdifferential is compact in finite dimension, by local boundedness and closedness of the graph of the Clarke subdifferential on the sampling ball. Therefore the minimum in \eqref{eq:ggsp_residual} is attained.

\paragraph{Incomparability at the same radius and accuracy.} Both directions can fail even when $\cX=[0,1]$. First take $f(t)=t$ and $0<x\le\delta<1/2$. Since $0$ is a feasible sample and $1-1=0$ with $-1\in N_{\cX}(0)$, we have $S_\delta(x)=0$. Here $M_{\gamma,\delta}(x)=\min\{1,x/\gamma\}>0$. Second, take $f(t)=-|t+a|$, $x=0$, and $0<a<\delta<1/2$. The ambient ball includes the kink at $-a$, and hence $0\in\partial_\delta f(0)$ and $M_{\gamma,\delta}(0)=0$. On the feasible sampling interval, the objective subgradient is $-1$ and all normals are nonpositive; hence $S_\delta(0)=1$. The first example uses the normal at a neighboring feasible point. The second uses objective information outside the feasible set.

\begin{proposition}[Relations between the two stationarity measures]\label{prop:ggsp_comparison}
Let $\cX$ be nonempty, closed, and convex, let $\bx\in\cX$, and fix $\gamma,\delta,\epsilon>0$. Put $R_\gamma=\max\{\delta,\gamma\epsilon\}$.
\begin{enumerate}[(i)]
\item If $N_\cX(\by)\subseteq N_\cX(\bx)$ for every $\by\in\ball_\delta(\bx)\cap\cX$, then $M_{\gamma,\delta}(\bx)\le S_\delta(\bx)$.
In particular, the normal-cone inclusion holds if $\ball_\delta(\bx)\subseteq\operatorname{int}(\cX)$, since the normal cone at every point in the ball is $\{\mathbf{0}\}$.
\item If $\ball_{r_{\mathrm{in}}}(\bc)\subseteq\cX$, $\delta<r_{\mathrm{in}}$, and $f$ is globally $L_f$-Lipschitz, then
\[
S_\delta(\bx)\le\epsilon\quad\Longrightarrow\quad
M_{\gamma,\delta}(\bx)\le\epsilon+
\sqrt{\frac{\delta\|\bx-\bc\|(L_f+\epsilon)}{\gamma(r_{\mathrm{in}}-\delta)}}.
\]
In particular, for a box centered at $\bc$ with half-widths $a_j$, one may take $r_{\mathrm{in}}=\min_j a_j$, provided that $\delta<\min_j a_j$.
\item If there exists
\[
\bv\in\operatorname{conv}\bigcup_{\by\in\ball_\delta(\bx)\cap\cX}\partial f(\by)
\qquad\text{with}\qquad
\|\bw_{\bv}\|\le\epsilon,
\]
then $S_{R_\gamma}(\bx)\le\epsilon$. In particular, if $\ball_\delta(\bx)\subseteq\cX$ and $\gamma\epsilon\le\delta$, then $M_{\gamma,\delta}(\bx)\le\epsilon$ implies $S_\delta(\bx)\le\epsilon$.
\item If $f$ is continuously differentiable on a neighborhood of $\ball_{R_\gamma}(\bx)$ and its gradient is $L_{\mathrm s}$-Lipschitz on that ball, then
\[
M_{\gamma,\delta}(\bx)\le\epsilon\quad\Longrightarrow\quad
S_{\gamma\epsilon}(\bx)\le\epsilon+L_{\mathrm s}(\delta+\gamma\epsilon).
\]
\end{enumerate}
\end{proposition}

The proof is given in Appendix~\ref{app:ggsp_comparison}. When $\ball_\delta(\bx)\subseteq\operatorname{int}(\cX)$ and $\gamma\epsilon\le\delta$, parts (i) and (iii) give equivalence at the same radius and accuracy.

Buzaglo and Kornowski \cite[Section 2 and Lemma 8.2]{BK26} use the residual
\[
\begin{aligned}
 \mathcal R_\delta^{f,\cX}(\bx)&:=\operatorname{dist}(\mathbf0,\partial_\delta f(\bx)+N_\cX(\bx))\\
 &=\min_{\bv\in\partial_\delta f(\bx)}
   \sup_{\gamma'>0}
   \frac{\|\bx-\proj_\cX(\bx-\gamma'\bv)\|}{\gamma'}
 \ge M_{\gamma,\delta}(\bx).
\end{aligned}
\]
For fixed $f$ and $\cX$, we abbreviate $\mathcal R_\delta^{f,\cX}$ to $\mathcal R_\delta$. The residual $\mathcal R_\delta$ retains the normal cone at $\bx$ while averaging objective subgradients over the ambient sampling ball. The measure $S_\delta$ applies the neighborhood convexification to the sums $\partial f(\by)+N_\cX(\by)$ at feasible samples. In the two examples above, $\mathcal R_\delta$ equals $1$ and $0$, respectively. Thus neither residual bounds the other in general. The supremum over projection steps gives a stronger requirement than the fixed-step criterion of Liu et al.

\begin{proposition}[Comparison of $\mathcal R_\delta$ and $S_\delta$]\label{prop:residual_comparison}
Let $\cX$ be nonempty, closed, and convex, let $f:\RR^d\to\RR$ be locally Lipschitz, and fix $\bx\in\cX$ and $\delta>0$.
\begin{enumerate}[(i)]
\item If $N_\cX(\by)\subseteq N_\cX(\bx)$ for every $\by\in\ball_\delta(\bx)\cap\cX$, then $\mathcal R_\delta(\bx)\le S_\delta(\bx)$.
\item If $\ball_\delta(\bx)\subseteq\cX$, then $S_\delta(\bx)\le\mathcal R_\delta(\bx)$. If $\ball_\delta(\bx)\subseteq\operatorname{int}(\cX)$, both residuals equal $\dist(\mathbf0,\partial_\delta f(\bx))$.
\item If $f$ is continuously differentiable on a neighborhood of $\ball_\delta(\bx)$ and its gradient is $L_{\mathrm s}$-Lipschitz on that ball, then
\[
S_\delta(\bx)\le\mathcal R_\delta(\bx)+L_{\mathrm s}\delta.
\]
\end{enumerate}
\end{proposition}
The proof is given in Appendix~\ref{app:ggsp_comparison}.

\paragraph{Goldstein KKT stationarity.}
Consider a problem with one explicit inequality constraint:
\begin{equation}\label{eq:primal_inequality}
 \min_{\bx\in\RR^d}f(\bx)\quad\text{s.t.}\quad g(\bx)\le0,
\end{equation}
where neither the objective $f$ nor the constraint $g$ is required to be smooth or convex. Let $\mathcal S:=\{\bx\mid g(\bx)\le0\}$ denote the induced feasible set. The constraint function $g$ is fixed independently of the objective construction.

The following definition of Goldstein KKT stationarity is used in the recent literature on inequality constraints \cite{GSL}. 

\begin{definition}[$(\delta,\epsilon,\eta)$-Goldstein KKT Stationarity \cite{GSL}]\label{def:standard_gkkt}
A feasible point $\bx \in \mathcal{S}$ is a $(\delta,\epsilon,\eta)$-Goldstein KKT (GKKT) stationary point for problem \eqref{eq:primal_inequality} if there exists a multiplier $\lambda \ge 0$ such that:
\[
    \dist\Big(\mathbf{0}, \, \partial_\delta f(\bx) + \lambda \partial_\delta g(\bx)\Big) \le \epsilon,
\]
and the relaxed complementary slackness condition holds:
\[
    \max_{\by \in \ball_\delta(\bx)} |\lambda g(\by)| \le \eta.
\]
\end{definition}

Unlike the set-based condition, this definition uses separate ambient Goldstein subdifferentials and a scalar multiplier. Its complementary condition is imposed on the full ball, not just on feasible points. The effect of neighborhood convexification must be considered together with this complementary-slackness requirement. For example, for $f(x)=-x$ and $g(x)=x$, the feasible point $x=0$ is an exact KKT point with multiplier $1$. For $x=0$ to be $(\delta,\epsilon,\eta)$-GKKT stationary with $0\le\epsilon<1$, a multiplier would have to satisfy both $\lambda\ge1-\epsilon$ and $\lambda\delta\le\eta$. No such multiplier exists if $\eta<(1-\epsilon)\delta$. Consequently, no general ordering with pointwise approximate KKT conditions follows merely from the use of Goldstein subdifferentials; the neighborhood radius and the tolerance for complementary slackness must also be compared.

\subsection{Geometry of constraint sets}
\label{subsec:convex_geo}

For the geometric lower bound in Section~\ref{sec:framework1}, assume that the closed feasible set is \textit{star-shaped with respect to a ball}. For a closed set $\cX\subseteq\RR^d$, its \textit{kernel}, denoted by $\ker\cX$ as in \cite[p.~703]{Breen80}, is
\[
 \ker\cX:=\{\bp\in\cX\mid[\bp,\by]\subseteq\cX\ \text{for every }\by\in\cX\}.
\]
Assume that, for some fixed center $\bc\in\cX$ and radius $r>0$,
\begin{equation}\label{eq:ball_kernel}
 \ball_r(\bc)\subseteq\ker\cX.
\end{equation}
Condition \eqref{eq:ball_kernel} means that $\cX$ is star-shaped with respect to every point of $\ball_r(\bc)$; equivalently, $\operatorname{conv}(\ball_r(\bc)\cup\{\by\})\subseteq\cX$ for every $\by\in\cX$. This is the closed-set version of the standard notion of a domain star-shaped with respect to a ball; see \cite[Section 2]{GS21}, \cite[Definition 1.2(ii)]{GGS15}, and \cite[Section 2]{CM10}. Every closed convex set containing $\ball_r(\bc)$ satisfies this condition. Specifically, it includes all convex bodies with nonempty interior \cite[Section 1.1]{schneider2014convex}, as well as unbounded full-dimensional convex sets and some nonconvex star-shaped sets.

For Lipschitz domains, the relation between this form of star-shapedness and outward-normal estimates is also noted in \cite[Remark 1.3]{GGS15}. The following proposition gives a version for limiting normals of a closed set.

\begin{proposition}[Normal-Cone Estimate]\label{prop:ball_kernel_normal}
If $\cX$ satisfies \eqref{eq:ball_kernel}, then
\[
 \langle\bn,\by-\bc\rangle\ge r\|\bn\|
 \qquad\text{for all }\by\in\cX,\ \bn\in N_\cX(\by).
\]
Consequently, if $\by\in\ball_\delta(\bx)\cap\cX$ and $\delta<r$, then
\[
 \langle\bn,\bx-\bc\rangle\ge(r-\delta)\|\bn\|\ge0.
\]
\end{proposition}
\begin{proof}
First, let $\bn\in\hat N_\cX(\by)$. For any $\bp\in\ball_r(\bc)$, the segment $\by+t(\bp-\by)$, $0\le t\le1$, lies in $\cX$. The defining regular-normal inequality, divided by $t$ and followed by $t\downarrow0$, gives $\langle\bn,\bp-\by\rangle\le0$. Choosing $\bp=\bc+r\bn/\|\bn\|$ when $\bn\ne0$ proves the first estimate for regular normals. The estimate passes to limiting normals by their defining sequences. The second estimate follows from $\|\bx-\by\|\le\delta$ and Cauchy--Schwarz.
\end{proof}

\subsection{Deterministic algorithms and oracle complexity}\label{subsec:oracle_model}

Following \cite{nogas,DY26}, at a query $\bq\in\RR^d$, a local oracle returns a locally Lipschitz function $O_f(\bq):\RR^d\to\RR$ satisfying
\[
 O_f(\bq)(\by)=f(\by)\qquad\text{for all }\by\in B(\bq,\nu)
\]
for some $\nu>0$, whose value is not supplied. The algorithm's internal computation, queries, stopping decision, and output depend only on the fixed input and the local information in preceding responses. Replacing responses by functions that agree with them near the corresponding query points cannot change the subsequent computation. Queries are unrestricted and may be repeated, and the returned point need not have been queried. All objective-oracle calls are counted; arbitrary finite internal computation is allowed between calls.

The constraint set or function $g$ is fixed and may be given in its entirety. The input includes $d$, a reference feasible point $\bc$, bounds $L_f,\Delta_f>0$, and accuracies $\delta,\epsilon>0$. Define
\[
 \mathcal F_\cX(L_f,\Delta_f;\bc)
 :=\left\{f:\RR^d\to\RR\ \middle|\ 
 \begin{array}{l}
  f\text{ is globally }L_f\text{-Lipschitz},\\
  \inf_{\bx\in\cX}f(\bx)>-\infty,\\
  f(\bc)-\inf_{\bx\in\cX}f(\bx)\le\Delta_f
 \end{array}\right\}.
\]
For a deterministic algorithm $A$, let $\tau_A(f;\bc,\delta,\epsilon)$ be the total number of calls if it returns a feasible $(\delta,\epsilon)$-Goldstein stationary point in finite time, and $+\infty$ otherwise. Following \cite{DY26}, the optimal worst-case oracle complexity is
\[
 \mathcal T_\cX(d,L_f,\Delta_f,\delta,\epsilon;\bc)
 :=\inf_A\ \sup_{f\in\mathcal F_\cX(L_f,\Delta_f;\bc)}
       \tau_A(f;\bc,\delta,\epsilon).
\]
The infimum is over all deterministic algorithms in this local-oracle model. Similarly, for a fixed constraint function $g$ and $\eta\ge0$, set $\mathcal S=\{g\le0\}$ and define $\tau_A(f;\bc,\delta,\epsilon,\eta)$ analogously, using Definition~\ref{def:standard_gkkt} as the valid-output criterion. The corresponding objective-oracle complexity is
\[
 \mathcal T_g(d,L_f,\Delta_f,\delta,\epsilon,\eta;\bc)
 :=\inf_A\ \sup_{f\in\mathcal F_{\mathcal S}(L_f,\Delta_f;\bc)}
       \tau_A(f;\bc,\delta,\epsilon,\eta).
\] Fixed parameters are suppressed below. A uniform $T$-call guarantee requires the same algorithm to return a valid point within $T$ calls for every admissible objective. The local-oracle lower bounds also hold for a value--subgradient oracle returning $f(\bq)$ and one Clarke subgradient in $\partial f(\bq)$.

\section{Complexity lower bound under set constraints}
\label{sec:framework1}

The normal vectors must be controlled without restricting where the algorithm can query. The construction combines a scaled family of linear functions with mutually orthogonal gradients and a resisting-oracle argument, following the deterministic lower-bound constructions in \cite{nemirovski1983problem,KS22,nogas,DNN}. The affine components are active only in interior neighborhoods, where the normals vanish. Near the boundary, the radial gradients have positive projections in the direction from the reference point to the neighborhood center, and adding normal vectors does not decrease these projections by Proposition~\ref{prop:ball_kernel_normal}. The feasible set is fixed and satisfies \eqref{eq:ball_kernel}; convexity and compactness are not otherwise required.

The precision-dependent bound below supplements the exponential bounds in Section~\ref{sec:exponential} by allowing the broader geometric condition \eqref{eq:ball_kernel}. Appendix~\ref{app:precision_proofs} gives the construction and the estimate of its initial objective-value gap.

\subsection{Oracle complexity lower bound and consequences}

\begin{theorem}[Precision-Dependent Oracle Complexity]\label{thm:set_complexity}
Let $\cX$ satisfy \eqref{eq:ball_kernel}, fix $L_f,\Delta_f>0$, and set
\[
 \ell:=\min\left\{L_f,\frac{2\Delta_f}{r}\right\}.
\]
For every integer $T\ge1$ satisfying
\[
 T+2\le d,\qquad
 0<\delta\le\frac{r}{16\sqrt{T+2}},\qquad
 0<\epsilon\le\frac{\ell}{16\sqrt{T+2}},
\]
and every deterministic algorithm, there is a globally $L_f$-Lipschitz convex objective with $f(\bc)-\inf_{\bx\in\cX}f(\bx)\le\Delta_f$ on which the algorithm cannot return a $(\delta,\epsilon)$-Goldstein stationary point within $T$ calls. The constructed function has a global minimizer in the interior of the feasible set.
\end{theorem}
\begin{proof}
Use Lemma~\ref{lem:general_convex_lower_bound_sqrtT} with $L=\ell$. Estimate \eqref{eq:initial_gap_upper} gives the prescribed objective-gap bound. The residual exceeds $\epsilon$ by \eqref{eq:set_residual_bound} at every feasible point returned after at most $T$ calls; an infeasible output is invalid. If the algorithm does not return after at most $T$ calls, the guarantee also fails.
\end{proof}

\begin{corollary}[Worst-Case Oracle Complexity]\label{cor:set_uniform}
Put
\[
 Q:=\min\left\{d,\frac{\ell^2}{16^2\epsilon^2},
                    \frac{r^2}{16^2\delta^2}\right\}.
\]
If $Q\ge3$, the optimal worst-case oracle complexity $\mathcal T_\cX$ is at least $\lfloor Q\rfloor-1$. Consequently,
\[
 \mathcal T_\cX
 =\Omega\left(\min\left\{d,\frac{\ell^2}{\epsilon^2},\frac{r^2}{\delta^2}\right\}\right).
\]
\end{corollary}
\begin{proof}
Take $T=\lfloor Q\rfloor-2$ in Theorem~\ref{thm:set_complexity}.
\end{proof}

\begin{remark}[Convex Domains and the Neighborhood Radius]\label{rem:set_scope}
Every full-dimensional closed convex set containing the reference ball satisfies \eqref{eq:ball_kernel}. Thus the theorem covers balls, boxes, simplices, and other convex bodies with nonempty interior, as well as some nonconvex star-shaped sets. Similar dimension conditions are used in \cite[Theorem 1 and Appendix B.1]{DNN} and \cite[Theorem 1]{nogas}, with $d\ge T+2$ and $d\ge T+1$, respectively. In our construction, one orthogonal gradient direction is assigned to each of the at most $T+2$ retained centers: the reference point, at most $T$ queries, and a possible returned point.
\end{remark}

The $\Omega(\epsilon^{-2})$ order on a ball already follows from classical convex minimization lower bounds for deterministic local oracles \cite[Section 4.4.3, Proposition 1]{nemirovski1983problem}; Proposition~\ref{prop:objective_gap} in Appendix~\ref{app:classical_precision_comparison} derives an objective-value gap bound from approximate stationarity. Kornowski \cite[Theorem 3.1]{KFO26} and Buzaglo and Kornowski \cite[Corollary 8.23]{BK26} establish $\Omega(\delta^{-1}\epsilon^{-3})$ first-order lower bounds under the respective accuracy, oracle, and algorithm assumptions discussed in Section~\ref{subsec:related_work}. Scaling their initial gap also gives $\epsilon^{-2}$ consequences in suitable parameter ranges. Theorem~\ref{thm:set_complexity} does not improve this precision order. It gives a uniform construction for every prescribed set satisfying \eqref{eq:ball_kernel}, including nonconvex and unbounded sets, and controls normal vectors throughout the sampling neighborhood. The oracle-complexity lower bound in Corollary~\ref{cor:set_uniform} depends jointly on $d,\epsilon,\delta$ and does not involve an outer-radius parameter. The convex constructed function has a global minimizer in the interior of the feasible set.

\section{Exponential oracle complexity lower bounds}\label{sec:exponential}

The main result is an exponential lower bound for set-constrained Goldstein stationarity on a class of convex product domains. The proof combines the Brouwer query lower bound with a finite-set interpolation that accounts for nearby normal vectors.

\subsection{Constraint geometry and lower bounds}

Assume an orthogonal block decomposition $\RR^d=\prod_{j=1}^k\RR^{d_j}$, with $d_j\ge1$. For the exponential lower bounds, let $\cX=\prod_{j=1}^k \cX_j$, where each $\cX_j$ is a convex body satisfying
\begin{equation}\label{eq:product_geometry}
\ball_{r_j}(\bc_j)\subseteq \cX_j\subseteq \ball_{R_j}(\bc_j),
\qquad r_j\ge r_{\mathrm{in}}>0,\qquad \frac{R_j}{r_j}\le K,
\qquad K\ge1.
\end{equation}
Then $\bc=(\bc_1,\ldots,\bc_k)$ and $\ball_{r_{\mathrm{in}}}(\bc)\subseteq \cX$. Product domains are standard in block-coordinate optimization \cite[Section 1]{BRZ25}. Inner and outer balls and their radius ratio are also used to quantify the geometry of a convex body \cite{LPST26}. Here the uniform ratio bound is imposed separately on the blocks. It gives a lower bound on the measure of the common directions used in Lemma~\ref{lem:constraint_cones}. A Euclidean ball can be taken as one block with $K=1$. An axis-aligned box also has $K=1$ when each interval of half-width $a_j>0$ is one block with $r_j=R_j=a_j$. The conditions also allow products of balls and ellipsoids with uniformly bounded axis ratios.

The Brouwer query bound in Lemma~\ref{lem:brouwer_query} supplies absolute constants $H_{\mathrm B}\ge1$, $0<\varepsilon_{\mathrm B}\le1$, $b_{\mathrm B}>0$, and an integer $d_{\mathrm B}\ge2$. For $d\ge d_{\mathrm B}$, put $Q_{\mathrm B}(d)=\lfloor\exp(b_{\mathrm B}d)\rfloor$. Define
\[
\beta_{\mathrm B}=\frac{\varepsilon_{\mathrm B}}{128},\qquad
\kappa=\frac1{32K},\qquad
\theta_{\mathrm B}=\min\left\{\frac1{64},\frac{\beta_{\mathrm B}}{4(H_{\mathrm B}+2)}\right\},\qquad
\rho_{\mathrm B}=\theta_{\mathrm B} r_{\mathrm{in}}.
\]

For an integer $T\ge0$, consider the conditions
\begin{equation}\label{eq:exp_parameters}
\begin{gathered}
d\ge d_{\mathrm B},\qquad T+2\le Q_{\mathrm B}(d),\qquad
0<\epsilon\le \frac{\beta_{\mathrm B} L_f}{8},\\
(T+2)^2\left(\frac{32\epsilon}{L_f}\right)^{d-1}
<\left(\frac{\kappa}{2}\right)^d,\qquad
0<\delta\le\min\left\{\frac{\rho_{\mathrm B}\epsilon}{2L_f},
                       \frac{\Delta_f}{2L_f}\right\}.
\end{gathered}
\end{equation}

\begin{theorem}[Set-constrained Goldstein stationarity]\label{thm:exp_set}
Fix $\cX$ satisfying \eqref{eq:product_geometry} and parameters satisfying \eqref{eq:exp_parameters}. For every deterministic algorithm using a local oracle, there is a globally $L_f$-Lipschitz objective for which the algorithm fails to return a feasible point with $S_\delta(\bx)\le\epsilon$ within $T$ objective-oracle calls. The objective can be chosen to satisfy
\begin{equation}\label{eq:exp_objective}
f(\bc)=0,\qquad
\min_{\RR^d} f=-2L_f\delta,\qquad
\operatorname*{argmin}_{\RR^d}f\cap\operatorname{int}(\cX)\ne\varnothing.
\end{equation}
\end{theorem}

\begin{corollary}[Fixed-Accuracy Dimension-Exponential Regime]\label{cor:exp_fixed_accuracy}
Fix $K\ge1$ and $r_{\mathrm{in}},L_f,\Delta_f>0$. There are constants $\epsilon_0,c_0>0$ and an integer $d_0$ such that for all $d\ge d_0$ and $0<\epsilon\le\epsilon_0 L_f$, the set-constrained lower bound applies with $T=\lfloor\exp(c_0d)\rfloor$ whenever
\[
0<\delta\le\min\left\{\frac{\rho_{\mathrm B}\epsilon}{2L_f},\frac{\Delta_f}{2L_f}\right\}.
\]
The worst-case number of objective-oracle calls is therefore $\exp(\Omega(d))$ at fixed admissible accuracies.
\end{corollary}
The proof is given in Appendix~\ref{app:exponential_estimates}.

\subsection{Proof of the exponential lower bound}\label{subsec:exponential_proof}

\begin{lemma}[Brouwer query lower bound]\label{lem:brouwer_query}
There exist absolute constants $H_{\mathrm B}\ge1$, $\varepsilon_{\mathrm B}\in(0,1]$, $b_{\mathrm B}>0$, and an integer $d_{\mathrm B}\ge2$ such that for $d\ge d_{\mathrm B}$, any deterministic algorithm that returns a point of $\ball_1^d(\mathbf{0})$ on every $H_{\mathrm B}$-Lipschitz map $\mathbf{F}:\ball_1^d(\mathbf{0})\to \ball_1^d(\mathbf{0})$ using at most $Q_{\mathrm B}(d)$ value queries has an instance on which its output $\bz$ satisfies
\begin{equation}\label{eq:brouwer_query}
\|\bz-\mathbf{F}(\bz)\|>\varepsilon_{\mathrm B}.
\end{equation}
\end{lemma}

This follows from the communication lower bound for the Local Brouwer problem in \cite[Definitions 2.1 and 5.1 and Theorem 5.2]{GKP22}. The conversion to value queries on the unit ball and the dimension adjustment are given in Appendix~\ref{app:exponential_estimates}.

\paragraph{Slopes of the local oracle responses.} Let $\mathbf{P}=\proj_{\ball_1^d(\mathbf{0})}$, $a_{\mathrm B}=r_{\mathrm{in}}/8$, and $\ell_{\mathrm B}=L_f/64$. For a map $\mathbf{F}$ in Lemma~\ref{lem:brouwer_query}, write $\bz_{\bx}=(\bx-\bc)/a_{\mathrm B}$, $t_{\bx}=\max\{1,\|\bz_{\bx}\|\}$, and define
\begin{equation}\label{eq:brouwer_field}
\mathbf{G}(\bx)=\ell_{\mathrm B}\,\frac{\bz_{\bx}-\mathbf{F}(\mathbf{P}(\bz_{\bx}))}{t_{\bx}}.
\end{equation}
A single value of $\mathbf{F}$ determines $\mathbf{G}(\bx)$. Since $\bz_{\bx}/t_{\bx}=\mathbf{P}(\bz_{\bx})$, both terms in $\mathbf{G}(\bx)/\ell_{\mathrm B}$ have norm at most $1$. As functions of $\bz_{\bx}$, the projection is $1$-Lipschitz and the product $\mathbf{F}(\mathbf{P}(\bz_{\bx}))/t_{\bx}$ is $(H_{\mathrm B}+1)$-Lipschitz, using the bounds and the Lipschitz constants of its factors. Thus
\[
\|\mathbf{G}(\bx)\|\le\frac{L_f}{32},\qquad
\Lip(\mathbf{G})\le H_{\mathbf G}:=\frac{L_f(H_{\mathrm B}+2)}{8r_{\mathrm{in}}}.
\]

For $\bz=\bz_{\bq}$ and $t=t_{\bq}$, the bound $\|\mathbf{F}(\mathbf{P}(\bz))\|\le1$ gives $t-1\le\|\bz-\mathbf{F}(\mathbf{P}(\bz))\|$. Hence
\[
\|\mathbf{P}(\bz)-\mathbf{F}(\mathbf{P}(\bz))\|
\le\frac{\|\bz-\mathbf{F}(\mathbf{P}(\bz))\|+t-1}{t}
\le\frac{2\|\mathbf{G}(\bq)\|}{\ell_{\mathrm B}}.
\]
In particular, $\|\mathbf{G}(\bq)\|\le\beta_{\mathrm B} L_f=\ell_{\mathrm B}\varepsilon_{\mathrm B}/2$ implies that $\mathbf{P}(\bz_{\bq})$ has Brouwer residual at most $\varepsilon_{\mathrm B}$. The vector $\mathbf{G}(\bq)$ is a positive multiple of $\bq-\by_0$, where $\by_0=\bc+a_{\mathrm B}\mathbf{F}(\mathbf{P}(\bz_{\bq}))\in\ball_{a_{\mathrm B}}(\bc)$. This direction will be used to control the normal vectors near $\bq$.

For a Borel set $A\subseteq\mathbb S^{d-1}$, define the normalized spherical measure by
\[
\mathfrak m_d(A):=\mathbb P\left\{\frac{\bz}{\|\bz\|}\in A\right\},
\qquad \bz\sim\mathcal N(\mathbf0,I_d).
\]
Thus $\mathfrak m_d$ is the uniform probability measure on $\mathbb S^{d-1}$. The same notation applies in every positive dimension; in particular, $\mathfrak m_1$ assigns mass $1/2$ to each point of $\{-1,1\}$.

\begin{lemma}[Common directions for nearby normal vectors]\label{lem:constraint_cones}
\leavevmode\newline
Fix $\boldsymbol\xi=(\boldsymbol\xi_1,\ldots,\boldsymbol\xi_k)\in\cX$, and define
\begin{equation}\label{eq:nearby_normals}
C_{\cX}(\boldsymbol\xi)
=\left\{\bu\ \middle|\ \langle\bu,\bn\rangle\ge0
\ \ \forall\,\by\in\ball_{4\rho_{\mathrm B}}(\boldsymbol\xi)\cap\cX,\
\forall\,\bn\in N_{\cX}(\by)\right\}.
\end{equation}
This is a full-dimensional closed convex cone satisfying $\mathfrak m_d(C_{\cX}(\boldsymbol\xi)\cap\unit)\ge\kappa^d$. For every $\by_0\in\ball_{a_{\mathrm B}}(\bc)$, $\boldsymbol\xi-\by_0\in C_{\cX}(\boldsymbol\xi)$. If $\bu\in\operatorname{int}(C_{\cX}(\boldsymbol\xi))$ is a unit vector, then $\boldsymbol\xi-t\bu\in\operatorname{int}(\cX)$ for $0<t\le\rho_{\mathrm B}$.
\end{lemma}
\begin{proof}
The cone is closed and convex by definition. Write $\by=(\by_1,\ldots,\by_k)\in\ball_{4\rho_{\mathrm B}}(\boldsymbol\xi)\cap\cX$. The inner ball in block $j$ gives, for every nonzero $\bn_j\in N_{\cX_j}(\by_j)$,
\[
\langle\bn_j,\boldsymbol\xi_j-\bc_j\rangle
\ge(r_j-4\rho_{\mathrm B})\|\bn_j\|>0.
\]
In each block with such a normal, $\boldsymbol\xi_j\ne\bc_j$. Put $\bv_j=(\boldsymbol\xi_j-\bc_j)/\|\boldsymbol\xi_j-\bc_j\|$. Since $\|\boldsymbol\xi_j-\bc_j\|\le R_j$ and $4\rho_{\mathrm B}\le r_{\mathrm{in}}/16$,
\[
\left\langle\frac{\bn_j}{\|\bn_j\|},\bv_j\right\rangle
\ge\frac{r_j-4\rho_{\mathrm B}}{R_j}\ge\frac1{2K},\qquad
\left\langle\frac{\bn_j}{\|\bn_j\|},\bw\right\rangle
\ge\frac1{2K}-\|\bw-\bv_j\|\ge\frac1{4K}
\]
for every unit $\bw$ with $\|\bw-\bv_j\|\le1/(4K)$. Let $C_j$ be the cone generated by this cap, and let $C_j=\RR^{d_j}$ in blocks with no nonzero nearby normal. The product normal-cone formula and Appendix~\ref{app:exponential_estimates} give
\[
\prod_{j=1}^k C_j\subseteq C_{\cX}(\boldsymbol\xi),\qquad
\mathfrak m_d(C_{\cX}(\boldsymbol\xi)\cap\unit)
\ge\prod_{j=1}^k\kappa^{d_j}=\kappa^d>0.
\]
Thus $C_{\cX}(\boldsymbol\xi)$ is full-dimensional. For $\by_0\in\ball_{a_{\mathrm B}}(\bc)$, the inner ball gives $\langle\bn,\boldsymbol\xi-\by_0\rangle\ge(r_{\mathrm{in}}-4\rho_{\mathrm B}-a_{\mathrm B})\|\bn\|$ for the normals in \eqref{eq:nearby_normals}. Thus $\boldsymbol\xi-\by_0\in C_{\cX}(\boldsymbol\xi)$.

Finally, an interior vector $\bu$ has strictly positive inner product with each nonzero vector defining the halfspaces of $C_{\cX}(\boldsymbol\xi)$. Put $\bz=\boldsymbol\xi-t\bu$ with $0<t\le\rho_{\mathrm B}$. Suppose $\bz\notin\operatorname{int}(\cX)$, and set $\by=\proj_{\cX}(\bz)$. Then $\by\in\partial\cX$ and $\|\by-\boldsymbol\xi\|\le t$. If $\bz\ne\by$, take $\bn=\bz-\by$. Otherwise take any nonzero supporting normal at $\by$. In either case $\bn\in N_{\cX}(\by)\setminus\{\mathbf0\}$ and
\[
-t\langle\bn,\bu\rangle
=\langle\bn,\bz-\by\rangle+\langle\bn,\by-\boldsymbol\xi\rangle\ge0,
\]
contrary to $\langle\bn,\bu\rangle>0$.
\end{proof}

\begin{lemma}[A direction compatible with the interpolation centers]\label{lem:history_direction}
Suppose \eqref{eq:exp_parameters} holds. Let $E=\{\bq_1,\ldots,\bq_N\}$ consist of distinct points with $N\le T+2$, and let $\bq\in E\cap \cX$ satisfy $\|\mathbf{G}(\bq)\|>\beta_{\mathrm B} L_f$. Define $L_0=L_f/8$, $\epsilon_1=4\epsilon$, and $\alpha=\epsilon_1/L_0=32\epsilon/L_f$. For the cone $C=C_{\cX}(\bq)$ of Lemma~\ref{lem:constraint_cones}, there is a unit $\bu\in\operatorname{int}(C)$ such that
\begin{equation}\label{eq:history_direction}
\left|\left\langle \bu,\frac{\bq_i-\bq_j}{\|\bq_i-\bq_j\|}\right\rangle\right|
<\sqrt{1-\alpha^2}\quad(i\ne j),\qquad
\langle \bu,\mathbf{G}(\bq_i)\rangle\ge \epsilon_1\quad(\bq_i\in \ball_{2\rho_{\mathrm B}}(\bq)).
\end{equation}
\end{lemma}

\begin{proof}
Put $\bp=\mathbf{G}(\bq)/\|\mathbf{G}(\bq)\|\in C$. Appendix~\ref{app:exponential_estimates} gives $\mathfrak m_d(\{\bu\in C\cap\unit\mid\langle\bu,\bp\rangle\ge1/\sqrt2\})\ge(\kappa/2)^d$. The directions violating the first inequality in \eqref{eq:history_direction} lie in two spherical caps for each pair of centers. By the cap estimate in the same appendix and \eqref{eq:exp_parameters}, their union has measure at most
\[
\binom{N}{2}\,2\alpha^{d-1}
\le(T+2)^2\left(\frac{32\epsilon}{L_f}\right)^{d-1}
<(\kappa/2)^d.
\]
Since $\mathfrak m_d(\partial C\cap\unit)=0$, there is a unit $\bu\in\operatorname{int}(C)$ satisfying the first inequality in \eqref{eq:history_direction} and $\langle\bu,\bp\rangle\ge1/\sqrt2$. For $\bq_i\in\ball_{2\rho_{\mathrm B}}(\bq)$,
\[
\langle \bu,\mathbf{G}(\bq_i)\rangle
\ge\frac{\beta_{\mathrm B} L_f}{\sqrt2}-2H_{\mathbf G}\rho_{\mathrm B}
\ge\beta_{\mathrm B} L_f\left(\frac1{\sqrt2}-\frac1{16}\right)>\frac{\beta_{\mathrm B} L_f}{2}\ge4\epsilon,
\]
by the definitions of $\theta_{\mathrm B},\rho_{\mathrm B},H_{\mathbf G}$ and the accuracy condition in \eqref{eq:exp_parameters}.
\end{proof}

\begin{lemma}[Realizing the local oracle responses]\label{lem:exp_interpolation}
Under the assumptions of Lemma~\ref{lem:history_direction}, suppose $\bc\in E$ and let $\bu\in\operatorname{int}(C)$ be the unit direction furnished by that lemma. There exists an $L_f$-Lipschitz $f$ satisfying \eqref{eq:exp_objective} that agrees with $\bx\mapsto\langle \mathbf{G}(\bq_i),\bx-\bq_i\rangle$ on a neighborhood of every $\bq_i\in E$. Moreover,
\begin{equation}\label{eq:objective_projection}
\langle \bu,\bv\rangle\ge4\epsilon\qquad
\text{for all }\bv\in\partial_\delta f(\bq).
\end{equation}
\end{lemma}

\begin{proof}
The construction adapts the finite minimum of penalized affine functions in \cite[Lemma 3.2]{DY26} to slopes that differ across centers. Define
\[
\mathbf Q_{\bu}=\mathbf I-\bu\bu^\top,\qquad
b_0=\sqrt{L_0^2-\epsilon_1^2},\qquad
\psi(\bz)=\epsilon_1\langle\bu,\bz\rangle+b_0\|\mathbf Q_{\bu}\bz\|.
\]
For $\bz=\bq_i-\bq_j\ne\mathbf0$, the first inequality in \eqref{eq:history_direction} gives
\[
b_0\|\mathbf Q_{\bu}\bz\|>\epsilon_1|\langle\bu,\bz\rangle|,
\qquad \psi(\bq_i-\bq_j)>0.
\]
Choose $s\in(\rho_{\mathrm B},2\rho_{\mathrm B})$ with $s\ne\|\bq_i-\bq\|$ for all $i$, and put $\mathcal D=B(\bq,s)$. Prescribe the slopes
\[
\ba_i=\begin{cases}
\mathbf G(\bq_i),&\bq_i\in \mathcal D,\\
\epsilon_1\bu,&\bq_i\notin \mathcal D.
\end{cases}
\]
By \eqref{eq:exp_parameters}, \eqref{eq:history_direction}, and $\|\mathbf G(\bq_i)\|\le L_f/32$, all $\ba_i$ satisfy $\|\ba_i\|\le L_0$ and $\langle\bu,\ba_i\rangle\ge \epsilon_1$. If $N\ge2$, take
\[
0<\zeta<\min_{i\ne j}\psi(\bq_i-\bq_j).
\]
If $N=1$, take any $\zeta>0$. Define
\[
F_i(\bx)=\max\left\{
\langle\ba_i,\bx-\bq_i\rangle,\ \psi(\bx-\bq_i)-\zeta
\right\},\qquad
\phi_0(\bx)=\min_{1\le i\le N}F_i(\bx).
\]
At $\bq_i$, the two terms of $F_i$ are $0$ and $-\zeta$, whereas $F_j(\bq_i)\ge \psi(\bq_i-\bq_j)-\zeta>0$ for $j\ne i$. Continuity and finiteness therefore give an open ball $U_i$ about each center on which $\phi_0(\bx)=\langle\ba_i,\bx-\bq_i\rangle$. Shrink the $U_i$ to be disjoint, each contained in $\mathcal D$ or disjoint from $\overline{\mathcal D}$ according to its center.

The orthogonal decomposition gives
\[
\begin{aligned}
|\psi(\bx)-\psi(\by)|
&\le \epsilon_1|\langle\bu,\bx-\by\rangle|
+b_0\|\mathbf Q_{\bu}(\bx-\by)\|\\
&\le\sqrt{\epsilon_1^2+b_0^2}\,\|\bx-\by\|
=L_0\|\bx-\by\|.
\end{aligned}
\]
Each $F_i$ is the maximum of two $L_0$-Lipschitz functions, and their finite minimum $\phi_0$ is also $L_0$-Lipschitz. For $t\ge0$, the affine term of $F_i(\bx+t\bu)$ increases by $t\langle\ba_i,\bu\rangle\ge \epsilon_1t$, and its other term increases by $\epsilon_1t$. Taking the maximum and then the minimum gives
\begin{equation}\label{eq:interpolation_growth}
\phi_0(\bx+t\bu)\ge\phi_0(\bx)+\epsilon_1t
\qquad(\bx\in\RR^d,\ t\ge0).
\end{equation}
At every differentiability point, \eqref{eq:interpolation_growth} gives $\langle\bu,\nabla\phi_0(\bx)\rangle\ge \epsilon_1$. Passing to limits of gradients and taking convex combinations yields
\begin{equation}\label{eq:interpolation_projection}
\langle\bu,\bv\rangle\ge \epsilon_1
\qquad(\bx\in\RR^d,\ \bv\in\partial\phi_0(\bx)).
\end{equation}

To restore the slopes at centers outside $\mathcal D$, choose a closed ball $\ball_{s_i}(\bq_i)\subseteq U_i$ at each such center. Let $\chi:[0,\infty)\to[0,1]$ be the piecewise affine cutoff, specifically,
\[
\chi(t)=\begin{cases}
1,&0\le t\le\frac12,\\
2(1-t),&\frac12<t<1,\\
0,&t\ge1.
\end{cases}
\]
Define
\[
h_i(\bx)=\chi\left(\frac{\|\bx-\bq_i\|}{s_i}\right)
\langle \mathbf{G}(\bq_i)-\epsilon_1\bu,\bx-\bq_i\rangle,\qquad
\phi=\phi_0+\sum_{\bq_i\notin \mathcal D}h_i.
\]
The product estimate, $\Lip(\chi)=2$, and $\epsilon_1\le L_f/32$ give
\[
\Lip(h_i)\le3\|\mathbf G(\bq_i)-\epsilon_1\bu\|
\le\frac{3L_f}{16}.
\]
Split any line segment at the boundaries of the disjoint support balls. On each subsegment, at most one correction term is nonzero, and the sum is Lipschitz with constant $3L_f/16$. The sum vanishes at the ball boundaries, so adding the estimates along the segment gives
\[
\Lip\left(\sum_{\bq_i\notin\mathcal D}h_i\right)\le\frac{3L_f}{16},
\qquad
\Lip(\phi)\le\frac{L_f}{8}+\frac{3L_f}{16}=\frac{5L_f}{16}.
\]
Since the cutoff equals $1$ near $\bq_i$, $\phi$ has slope $\mathbf{G}(\bq_i)$ and value zero at each center, with affine agreement on a neighborhood. On $\mathcal D$, the correction terms vanish and $\phi=\phi_0$.

Set $f(\bx)=\max\{\phi(\bx),-2L_f\delta\}$. The truncation preserves the local responses because $\phi(\bq_i)=0$ and $\delta>0$. Let $\bx^*=\bq-(s/2)\bu$. Since $s/2<\rho_{\mathrm B}$, Lemma~\ref{lem:constraint_cones} puts $\bx^*$ in $\operatorname{int}(\cX)$. Moreover, \eqref{eq:interpolation_growth}, $\phi_0(\bq)=0$, and \eqref{eq:exp_parameters} give
\[
\phi(\bx^*)=\phi_0(\bx^*)\le -\frac{\epsilon_1s}{2}
<-2\epsilon\rho_{\mathrm B}<-2L_f\delta.
\]
Thus $f$ is constant at its global minimum on a neighborhood of $\bx^*$, and $f(\bc)=0$ because $\bc\in E$. Its initial objective-value gap is $2L_f\delta\le\Delta_f$ by \eqref{eq:exp_parameters}. The same condition gives $\delta<\rho_{\mathrm B}$ and hence $\ball_\delta(\bq)\subseteq \mathcal D$. On this ball,
\[
\phi(\by)\ge\phi(\bq)-L_f\|\by-\bq\|
\ge-L_f\delta>-2L_f\delta.
\]
By continuity, $f=\phi_0$ in a neighborhood of each sampling point $\by$. Applying \eqref{eq:interpolation_projection} at these points and taking convex combinations gives \eqref{eq:objective_projection}.
\end{proof}

\paragraph{Oracle simulation.} A deterministic optimization algorithm using at most $T$ calls would give a Brouwer algorithm using at most $T$ value queries. Suppose a deterministic algorithm returns a valid point after at most $T$ objective-oracle calls on every admissible objective. Fix its constraint information and the other input parameters. Given a map $\mathbf{F}$ as in Lemma~\ref{lem:brouwer_query} through a value oracle, simulate the algorithm by returning the affine function $\bx\mapsto\langle \mathbf{G}(\bq_i),\bx-\bq_i\rangle$ at each query $\bq_i$. Each call uses at most one value query to $\mathbf{F}$, and repeated queries return the same affine function.

For any $\mathbf{F}$ and any prefix of at most $T$ calls, fix $\bq_0\in\partial\cX$. Since $\|\bq_0-\bc\|\ge r_{\mathrm{in}}=8a_{\mathrm B}$, \eqref{eq:brouwer_field} gives $\|\mathbf{G}(\bq_0)\|\ge7\ell_{\mathrm B}/8>\beta_{\mathrm B} L_f$. Apply Lemmas~\ref{lem:history_direction} and~\ref{lem:exp_interpolation} with $\bq=\bq_0$ and with the distinct points among the queries in the prefix, $\bc$, and $\bq_0$ as centers. The resulting objective satisfies \eqref{eq:exp_objective} and agrees with the simulated responses near the query points. By the local-oracle model, the algorithm on this objective makes the same decisions from that prefix. If the simulation ever requests a $(T+1)$-st call, returns an infeasible point, or neither queries again nor returns after such a prefix, the algorithm fails on the corresponding objective. The assumed guarantee therefore forces a feasible output $\bq$ within $T$ calls for every $\mathbf{F}$.

The Brouwer algorithm outputs $\mathbf{P}((\bq-\bc)/a_{\mathrm B})$. This requires no further value query. Thus the number of queries to $\mathbf{F}$ is at most $T\le Q_{\mathrm B}(d)$. By Lemma~\ref{lem:brouwer_query}, there is an $\mathbf{F}$ for which the output has residual greater than $\varepsilon_{\mathrm B}$. The estimate following \eqref{eq:brouwer_field} gives $\|\mathbf{G}(\bq)\|>\beta_{\mathrm B} L_f$ on this instance. Apply Lemmas~\ref{lem:history_direction} and~\ref{lem:exp_interpolation} to the distinct points among the queries, $\bc$, and $\bq$. The last two points enter the construction of $f$ as interpolation centers after the simulated run. For the resulting fixed objective $f$, each simulated affine response is a valid representative $O_f(\bq_i)$, and the deterministic algorithm returns the same $\bq$.

\begin{proof}[Proof of Theorem~\ref{thm:exp_set}]
Use $C_{\cX}(\bq)$. For each $\by\in\ball_\delta(\bq)\cap\cX$, \eqref{eq:objective_projection} gives $\langle\bu,\bv\rangle\ge4\epsilon$ for every $\bv\in\partial f(\by)$. Since $\delta<4\rho_{\mathrm B}$, \eqref{eq:nearby_normals} gives $\langle\bu,\bn\rangle\ge0$ for every $\bn\in N_{\cX}(\by)$. Taking convex combinations as in Definition~\ref{def:constrained_goldstein} gives $S_\delta(\bq)\ge4\epsilon>\epsilon$, a contradiction.
\end{proof}

\subsection{The fixed-step criterion}

\begin{corollary}[Fixed-step generalized Goldstein stationarity]\label{cor:ggsp_exp}
Fix $\cX$ and parameters satisfying \eqref{eq:product_geometry} and \eqref{eq:exp_parameters}. Suppose the prescribed step $\gamma>0$ satisfies
\begin{equation}\label{eq:ggsp_step}
\gamma\epsilon\le4\rho_{\mathrm B}.
\end{equation}
For every deterministic algorithm using a local oracle, there exists a globally $L_f$-Lipschitz objective satisfying \eqref{eq:exp_objective} for which the algorithm fails to return a feasible point with $M_{\gamma,\delta}(\bx)\le\epsilon$ within $T$ objective-oracle calls. In particular, under the conditions of Corollary~\ref{cor:exp_fixed_accuracy}, the lower bound is $\exp(\Omega(d))$.
\end{corollary}

\begin{proof}
Use the oracle simulation in Section~\ref{subsec:exponential_proof} with $C_{\cX}$, taking the prescribed $\gamma$ as part of the fixed input and $M_{\gamma,\delta}(\bx)\le\epsilon$ as the valid-output condition. The consistency and stopping argument there applies because it only requires an admissible objective and a feasible output. Lemma~\ref{lem:exp_interpolation} therefore gives an objective on which the algorithm returns a feasible point $\bq$ satisfying \eqref{eq:objective_projection}, with $\bu\in C_{\cX}(\bq)$.

If $M_{\gamma,\delta}(\bq)\le\epsilon$, choose $\bv\in\partial_\delta f(\bq)$ with $\|\bw_{\bv}\|\le\epsilon$. Its projection point satisfies $\|\bq-\bp_{\bv}\|=\gamma\|\bw_{\bv}\|\le\gamma\epsilon\le4\rho_{\mathrm B}$. Thus $\bn_{\bv}$ in \eqref{eq:ggsp_projection} has nonnegative projection onto $\bu$ by \eqref{eq:nearby_normals}, and
\[
\|\bw_{\bv}\|\ge\langle \bu,\bw_{\bv}\rangle
=\langle \bu,\bv\rangle+\langle \bu,\bn_{\bv}\rangle
\ge4\epsilon,
\]
a contradiction. At the interior global minimizer constructed in Lemma~\ref{lem:exp_interpolation}, the function is locally constant, and $\bv=\mathbf0$ gives zero projection residual. Thus the target criterion has a feasible solution.
\end{proof}

Condition \eqref{eq:ggsp_step} ensures that any vector satisfying the target residual has its projection point within $4\rho_{\mathrm B}$ of the returned point, where \eqref{eq:nearby_normals} applies. On a compact feasible set, every feasible point satisfies the criterion if $\gamma\epsilon\ge\diam(\cX)$, since $M_{\gamma,\delta}(\bx)\le\diam(\cX)/\gamma$.

\section{Further results for explicit inequality constraints}\label{sec:inequality_constraints}

\subsection{Constraint functions}

For explicit inequality constraints, this section gives a precision-dependent lower bound under convexity.

\begin{assumption}[$\mathcal G$: Convex Coercive Slater Constraints]\label{assump:g}
A constraint function $g$ belongs to $\mathcal G$ if it is finite-valued and convex on $\RR^d$, is coercive, and has a strictly feasible reference point $\bc$:
\[
 g(\bx)\longrightarrow+\infty\quad\text{as }\|\bx\|\longrightarrow\infty,
 \qquad g(\bc)=-m<0.
\]
\end{assumption}

A finite convex function is locally Lipschitz \cite[Theorem 10.4]{R70}. Coercivity gives a compact feasible sublevel set, while strict feasibility gives nonempty interior. Hence there are $h,L_g>0$ such that $g$ is $L_g$-Lipschitz on $\mathcal S_h:=\{\bx\mid\dist(\bx,\mathcal S)<h\}$. Smoothness is not required; constraints such as $g(\bx)=\|\bx\|_1-1$ are admissible. Squared-norm ball constraints $g(\bx)=\|\bx-\bc\|^2-R_{\mathrm b}^2$, with $R_{\mathrm b}>0$, are also allowed, with the Lipschitz bound required only on $\mathcal S_h$.

Under the finite-valued convexity and Slater assumptions, coercivity is equivalent to compactness of $\mathcal S$. To see the converse implication, choose $D>\max_{\bx\in\mathcal S}\|\bx-\bc\|$ and put $m_D:=\min_{\|\by-\bc\|=D}g(\by)>0$. Convexity along the ray from $\bc$ to $\bx$ yields
\[
 g(\bx)\ge-m+\frac{\|\bx-\bc\|}{D}(m_D+m)
 \qquad\text{if }\|\bx-\bc\|>D.
\]
Thus compactness of the feasible sublevel set already implies coercivity.

\subsection{Precision-dependent lower bounds under convex constraints}

The functions
\[
\begin{aligned}
 &g_{\rm ball}(\bx)=\|\bx-\bc\|-R_{\mathrm b},\\
 &g_{\rm box}(\bx)=\max_{1\le i\le d}\{|x_i-c_i|-a_i\},\\
 &\quad R_{\mathrm b}>0,\ a_i>0\ (1\le i\le d),
\end{aligned}
\]
both belong to $\mathcal G$ and are globally $1$-Lipschitz.

\begin{theorem}[Oracle Complexity under Convex Functional Constraints]\label{thm:complexity_lb_inequality}
Let $g\in\mathcal G$, use the parameters $h,L_g,\bc,m$ above, and let
\[
 r:=\min\left\{\frac{h}{2},\frac{m}{2L_g}\right\},\qquad
 \ell:=\min\left\{L_f,\frac{2\Delta_f}{r}\right\}.
\]
For every integer $T\ge1$ with $T+2\le d$, if
\begin{equation}\label{eq:convex_functional_accuracy}
 0<\delta\le\frac{r}{16\sqrt{T+2}},\qquad
 0<\epsilon\le\frac{\ell}{16\sqrt{T+2}},\qquad
 0\le\eta\le\frac{r\epsilon}{8},
\end{equation}
then every deterministic algorithm has a globally $L_f$-Lipschitz convex objective with initial objective-value gap at most $\Delta_f$ on which it cannot return a $(\delta,\epsilon,\eta)$-GKKT stationary point within $T$ calls. The constructed function has a global minimizer in the interior of the feasible set.
\end{theorem}

The proof is given in Appendix~\ref{app:precision_proofs}. It uses convexity, the Slater point, and the stated local Lipschitz bounds.

\begin{corollary}[Convex Functional Precision Regime]\label{cor:convex_precision}
If $0\le\eta\le r\epsilon/8$ and
\[
 Q:=\min\left\{d,\frac{\ell^2}{16^2\epsilon^2},
                    \frac{r^2}{16^2\delta^2}\right\}\ge3,
\]
the optimal worst-case objective-oracle complexity is at least $\lfloor Q\rfloor-1$. When $d\ge\ell^2/(16^2\epsilon^2)$ and $\delta\le r\epsilon/\ell$, the required number of objective-oracle calls is $\Omega(\ell^2/\epsilon^2)$.
\end{corollary}
\begin{proof}
Take $T=\lfloor Q\rfloor-2$ in Theorem~\ref{thm:complexity_lb_inequality}.
\end{proof}

Theorem~\ref{thm:complexity_lb_inequality} extends the construction of Section~\ref{sec:framework1} to every prescribed $g\in\mathcal G$, using convexity, the Slater point, and complementary slackness. The $\epsilon^{-2}$ order, even with $\eta$ independent of $\delta$, already follows on a ball from the classical convex minimization lower bound; see Appendix~\ref{app:classical_precision_comparison}. The supplementary result treats any finite-valued convex coercive constraint with a Slater point, using local Lipschitz and strict-feasibility parameters instead of an additional outer-radius parameter. These parameters may depend on the dimension. The construction also preserves convexity of the objective and an interior global minimizer.

\section{Discussion and conclusion}
\label{sec:discussion}

The main result, Theorem~\ref{thm:exp_set} and Corollary~\ref{cor:exp_fixed_accuracy}, gives an $\exp(\Omega(d))$ lower bound for set-constrained Goldstein stationarity on convex product domains. Corollary~\ref{cor:ggsp_exp} gives the fixed-accuracy exponential lower bound for generalized Goldstein stationarity with a fixed projection parameter satisfying \eqref{eq:ggsp_step}.

The supplementary precision-dependent bounds in Corollaries~\ref{cor:set_uniform} and~\ref{cor:convex_precision} apply to broader set and functional constraints. Since the constructed functions have global minimizers in the interior of the feasible set, the lower bounds concern the information required to find stationary points rather than a failure of the stationarity condition to have a solution. The precision order on balls is classical; the supplementary contribution is the uniform construction over the broader prescribed constraints, with the stated parameter dependence.

\paragraph{Implications for algorithms.}
The model allows the algorithm to use the entire feasible set or constraint function. Exact projection, when available, or a change of the algorithm's internal step size or search rule therefore does not remove these oracle-complexity lower bounds while the algorithm remains deterministic and obtains objective information through the same local oracle. Here $\gamma$ is a prescribed parameter of the generalized Goldstein stationarity criterion. The precision-dependent lower bounds continue to hold under objective convexity. Examples of useful additional structure include sharpness for weakly convex objectives \cite{SMW}, error-bound conditions for convex objectives \cite{RSG}, and controlled smoothness or white-box structure in unconstrained problems \cite{DNN}.

Several questions remain open. For the set-constrained problem, it is natural to ask whether the exponential lower bound extends beyond convex product domains, and whether a sharper dependence on the stationarity accuracy can be established. For inequality constraints, it remains to establish exponential lower bounds under general assumptions on the constraint function.

\appendix
\section{Proofs for the precision-dependent lower bounds}\label{app:precision_proofs}

\subsection{Construction with mutually orthogonal gradients}

\begin{lemma}[Lower Bound on the Stationarity Measure]\label{lem:general_convex_lower_bound_sqrtT}
Let $\cX$ satisfy \eqref{eq:ball_kernel} and let $T\ge1$ be an integer with $T+2\le d$. For any $L>0$, fix an arbitrary deterministic algorithm and its input, and denote its $t$-th query point by $\bx_t$ whenever that query is made. There exists a globally $L$-Lipschitz convex function $f:\RR^d\to\RR$ such that every feasible queried point $\bx_t$ with $t\le T$ satisfies
\begin{equation}\label{eq:set_residual_bound}
 \dist(\mathbf0,\partial_\delta^\cX f(\bx_t))
 >\frac{L}{16\sqrt{T+2}},
 \qquad 0<\delta\le\delta_0:=\frac{r}{16\sqrt{T+2}}.
\end{equation}
The same bound holds at any feasible point returned within $T$ calls, whether or not it was queried. Moreover, $f$ has a global minimizer in $B(\bc,r/2)$.
\end{lemma}

\begin{proof}
\par\medskip\noindent\textbf{Step 1: Geometric Normalization and Parameter Selection.}\par\nobreak\smallskip
By translational invariance, take $\bc=\mathbf0$. Define
\[
 \mu:=\frac{L}{4},\qquad
 \zeta:=\frac{Lr}{64(T+2)^{\frac{3}{2}}},\qquad
 \delta_0:=\frac{r}{16\sqrt{T+2}},\qquad
 f_0(\bx):=L\left(\|\bx\|-\frac{r}{2}\right).
\]
Choose $\nu:=\zeta/(3L)$ for the oracle-consistency neighborhoods.

\par\medskip\noindent\textbf{Step 2: Recursive Objective Construction.}\par\nobreak\smallskip
Insert an auxiliary point $\bq_1=\mathbf0$, which is not an oracle call, choose a unit vector $\bv_1$, and set
\[
 f_1(\bx):=\mu\langle\bx,\bv_1\rangle-\zeta.
\]
After the first $j-1$ centers have been inserted, let
\[
 \mathcal V_{j-1}:=\operatorname{span}\{\bv_1,\ldots,\bv_{j-1}\}.
\]
If the algorithm requests its $(j-1)$-th query $\bx_{j-1}$, insert its location as $\bq_j:=\bx_{j-1}$. Obtain $\bv_j$ by projecting $\bq_j$ onto $\mathcal V_{j-1}^\perp$ and normalizing when the projection is nonzero; otherwise choose any unit vector in that orthogonal complement. Then
\[
 \bv_j\perp\mathcal V_{j-1},\qquad
 \langle\bq_j,\bv_j\rangle\ge0,\qquad
 \bq_j\in\mathcal V_j.
\]
Define
\[
 f_j(\bx):=\mu\langle\bx-\bq_j,\bv_j\rangle-j\zeta,\qquad
 F_j(\bx):=\max\{f_0(\bx),f_1(\bx),\ldots,f_j(\bx)\}.
\]
At a new query location, return $F_j$; at a repeated location, reuse its first response.

Retain the query history up to a return, a request for call $T+1$, or internal nontermination. If a point is returned within $T$ calls, insert it as one further center using the same projection procedure, whether or not it was queried. Including $\bq_1$, there are at most $T+2$ centers. Complete any remaining centers with $\bq_j=\mathbf0$ and unused orthogonal unit vectors, which exist since $T+2\le d$. These extra centers are not oracle calls; the nontermination case requires no detection procedure.

\par\medskip\noindent\textbf{Step 3: Final Objective and Oracle Consistency.}\par\nobreak\smallskip
The final global objective is
\[
 f(\bx):=\max\{f_0(\bx),f_1(\bx),\ldots,f_{T+2}(\bx)\}.
\]
Since $\Lip(f_j)=\mu\le L$ and $\Lip(f_0)=L$, $f$ is globally $L$-Lipschitz. It is convex and coercive.

Fix a center $\bq_i$ and a future index $j>i$. Orthogonality and $\bq_i\in\mathcal V_i$ give, for $\|\by-\bq_i\|<\nu$,
\[
 f_j(\by)
 =\mu\langle\by-\bq_i,\bv_j\rangle
   -\mu\langle\bq_j,\bv_j\rangle-j\zeta
 \le\mu\nu-j\zeta,
\]
whereas
\[
 F_i(\by)\ge f_i(\by)\ge-\mu\nu-i\zeta.
\]
Since $2\mu\nu<\zeta$, future components are strictly inactive on this neighborhood, so $f=F_i$ on $B(\bq_i,\nu)$. Thus $f$ supplies the same local information at every query, including repeated locations. Determinism and locality reproduce the retained execution, including any return within $T$ calls.

At the auxiliary center, $f(\mathbf0)=-\zeta<0$. If $\|\bx\|\ge r/2$, then $f(\bx)\ge f_0(\bx)\ge0$. Coercivity therefore gives a global minimizer in $B(\mathbf0,r/2)\subseteq\operatorname{int}\cX$.

\par\medskip\noindent\textbf{Step 4: Active Components and Radial Confinement.}\par\nobreak\smallskip
Fix any retained center $\bx=\bq_i$ and $0<\delta\le\delta_0$. For every $\by\in\ball_\delta(\bx)$,
\begin{equation}\label{eq:set_neighborhood_floor}
 f(\by)\ge f_i(\by)\ge-i\zeta-\mu\delta
 \ge-(T+2)\zeta-\mu\delta.
\end{equation}
If a linear component $f_k$ is active at $\by$, then
\[
 \begin{aligned}
 \mu\langle\bx,\bv_k\rangle
 &=f_k(\by)+\mu\langle\bq_k,\bv_k\rangle+k\zeta
   +\mu\langle\bx-\by,\bv_k\rangle\\
 &\ge-(T+2)\zeta-2\mu\delta.
 \end{aligned}
\]
Also $f_k(\by)\le\mu\|\by\|$. Activity implies
\[
 L\left(\|\by\|-\frac{r}{2}\right)\le\mu\|\by\|,
 \qquad\text{and thus}\qquad \|\by\|\le\frac{2r}{3}.
\]
If any linear component is active anywhere in $\ball_\delta(\bx)$, every point in that ball has norm at most $2r/3+2\delta<r$. The ball is then contained in $\operatorname{int}\cX$, and all feasible-set normals there vanish.

If the radial component is active at $\by$, \eqref{eq:set_neighborhood_floor} gives
\[
 \|\by\|\ge\frac{r}{2}-\frac{(T+2)\zeta}{L}-\frac{\mu\delta}{L}
 \ge\frac r2-\frac{r}{32\sqrt{T+2}}>0.
\]
Its gradient is therefore $\bg_0(\by)=L\by/\|\by\|$, and
\[
 \langle\bg_0(\by),\bx\rangle
 \ge L(\|\by\|-\delta)
 \ge\frac{Lr}{2}-(T+2)\zeta-\frac{5L\delta}{4}.
\]

\par\medskip\noindent\textbf{Step 5: Convex Hull Aggregation and Zone Partition.}\par\nobreak\smallskip
By the definition of the convex hull, any $\br\in\partial_\delta^\cX f(\bx)$ can be written as a finite convex combination
\[
 \br=\sum_{j=1}^J\theta_j(\bg_j+\bz_j),\qquad
 \theta_j\ge0,\quad\sum_j\theta_j=1,
\]
where $\by_j\in\ball_\delta(\bx)\cap\cX$, $\bg_j\in\partial f(\by_j)$, and $\bz_j\in N_\cX(\by_j)$. Proposition~\ref{prop:ball_kernel_normal} gives $\langle\bz_j,\bx\rangle\ge0$.

The maximum-subdifferential rule for convex functions gives
\begin{equation}\label{eq:bgj-decomp}
 \bg_j=\beta_j\bg_0(\by_j)
       +(1-\beta_j)\mu\sum_{k=1}^{T+2}\lambda_{jk}\bv_k,
\end{equation}
with $\beta_j\in[0,1]$, $\lambda_{jk}\ge0$, $\sum_k\lambda_{jk}=1$, and $(1-\beta_j)\lambda_{jk}=0$ whenever $f_k$ is inactive. A radial term with zero coefficient is omitted. Define
\begin{equation}\label{eq:radial_weights}
 \beta:=\sum_j\theta_j\beta_j,\qquad
 w_k:=\sum_j\theta_j(1-\beta_j)\lambda_{jk},\qquad
 \tilde{\bg}_0:=\sum_j\theta_j\beta_j\bg_0(\by_j).
\end{equation}
Then
\[
 \sum_{k=1}^{T+2}w_k=1-\beta,\qquad \|\tilde{\bg}_0\|\le\beta L.
\]
Writing $\bar{\bg}:=\sum_j\theta_j\bg_j$ gives
\begin{equation}\label{eq:bgbar-decomp}
 \bar{\bg}=\tilde{\bg}_0+\mu\sum_{k=1}^{T+2}w_k\bv_k.
\end{equation}
Split the argument at $\beta^*:=1/(4\sqrt{T+2})$.

\par\medskip\noindent\textbf{Zone A ($\beta\le\beta^*$): Linear-Component Dominance.}\par\nobreak\smallskip
Since $\beta<1$, a linear component is active at some sampling point. Step 4 shows that the whole neighborhood lies in the interior; therefore $\br=\bar{\bg}$. Let $I:=\{k\in\{1,\ldots,T+2\}\mid w_k>0\}$ and set
\[
 \bu:=\frac{\sum_{k\in I}\bv_k}{\sqrt{\#I}},\qquad \|\bu\|=1.
\]
For each $k\in I$, $f_k$ is active at some sampling point. Step~4 gives, whenever the radial component is active at $\by$,
\[
 \begin{aligned}
  \langle\bu,\by\rangle
  &\ge-\sqrt{\#I}\left(\frac{(T+2)\zeta}{\mu}+2\delta\right)-\delta\\
  &\ge-\frac{3r}{16}-\frac{r}{16\sqrt{T+2}}.
 \end{aligned}
\]
The norm estimate for such $\by$ and
\[
 \frac12\left(\frac r2-\frac{r}{32\sqrt{T+2}}\right)
 -\left(\frac{3r}{16}+\frac{r}{16\sqrt{T+2}}\right)
 =\frac r{16}-\frac{5r}{64\sqrt{T+2}}>0
\]
imply $\langle\bu,\bg_0(\by)\rangle>-L/2$. Since $\beta\le1/(4\sqrt{T+2})<1/6$, orthogonality gives
\begin{equation}\label{eq:zone_a_bound}
 \begin{aligned}
  \|\br\|&\ge\langle\bu,\bar{\bg}\rangle\\
  &\ge\frac{L(1-\beta)}{4\sqrt{T+2}}-\frac{\beta L}{2}
  \ge\frac{L}{12\sqrt{T+2}}.
 \end{aligned}
\end{equation}

\par\medskip\noindent\textbf{Zone B ($\beta>\beta^*$): Radial Activation.}\par\nobreak\smallskip
First, consider the case $\|\bx\|\le r$. The pointwise estimates in Step 4 give
\begin{equation}\label{eq:zone_b_radial}
 \begin{aligned}
  \langle\bar{\bg},\bx\rangle
  &\ge\frac{\beta Lr}{2}-(T+2)\zeta
     -\left(\frac12+\frac{3\beta}{4}\right)L\delta\\
  &>\frac{9Lr}{128\sqrt{T+2}}.
 \end{aligned}
\end{equation}
The preceding expression is increasing in $\beta$. At $\beta=\beta^*$, it is at least
\[
 \left(\frac5{64}-\frac3{256\sqrt{T+2}}\right)
 \frac{Lr}{\sqrt{T+2}}>\frac{9Lr}{128\sqrt{T+2}},
\]
since $\sqrt{T+2}>3/2$. Adding normals does not reduce this projection. Hence
\[
 \|\br\|\ge\frac{\langle\br,\bx\rangle}{r}
 >\frac{9L}{128\sqrt{T+2}}.
\]
If $\bx=\mathbf0$, \eqref{eq:zone_b_radial} contradicts $\langle\bar{\bg},\mathbf0\rangle=0$.

If $\|\bx\|>r$, then every $\by\in\ball_\delta(\bx)$ has norm greater than $r-\delta>2r/3$. Only the radial component is active. With $\bu_{\bx}:=\bx/\|\bx\|$,
\begin{equation}\label{eq:far_radial}
 \langle\bar{\bg},\bu_{\bx}\rangle
 \ge L\frac{\|\bx\|-\delta}{\|\bx\|+\delta}
 >\frac{L}{2}.
\end{equation}
The last inequality follows from $\delta<r/3<\|\bx\|/3$. The normal projection is nonnegative; therefore $\|\br\|>L/2$. Since each queried point $\bx_t$ with $t\le T$ and any point returned within $T$ calls is among the retained centers, taking the infimum over $\br\in\partial_\delta^\cX f(\bx)$ gives
\[
 \dist(\mathbf0,\partial_\delta^\cX f(\bx))
 \ge\frac{9L}{128\sqrt{T+2}}>\frac{L}{16\sqrt{T+2}}.
\]
This proves \eqref{eq:set_residual_bound} and the claim for returned points.

The objective estimates \eqref{eq:zone_a_bound}, \eqref{eq:zone_b_radial}, and \eqref{eq:far_radial} do not require the objective sampling points to be feasible: they hold for $\bar{\bg}\in\partial_\delta f(\bx)$ with the same active-weight decomposition. Feasibility was used only when adding the normal vectors. This observation will be used for functional constraints.
\end{proof}

\subsection{Estimate of the function-value gap}

For the prescribed reference point $\bc$, put
\[
 \Delta_0:=f(\bc)-\inf_{\bx\in\cX}f(\bx).
\]

Again, take $\bc=\mathbf0$. Since $f_0(\mathbf0)=-Lr/2$, the first linear component satisfies $f_1(\mathbf0)=-\zeta$, and every subsequent component satisfies $f_k(\mathbf0)\le-k\zeta$, it follows that $f(\mathbf0)=-\zeta$. Also $f(\bx)\ge f_0(\bx)\ge-Lr/2$ for every $\bx$. Therefore
\begin{equation}\label{eq:initial_gap_upper}
 \Delta_0\le\frac{Lr}{2}-\zeta<\frac{Lr}{2}.
\end{equation}
The infimum over $\cX$ equals the ambient infimum since the global minimizer lies in $B(\mathbf0,r/2)$.

\subsection{Convex functional constraints}

\begin{proof}[Proof of Theorem~\ref{thm:complexity_lb_inequality}]
Without loss of generality, translate the reference point to the origin, so that $g(\mathbf0)=-m$. The reference ball is strictly feasible and lies in the kernel of the convex set $\mathcal S$. Use the construction of Lemma~\ref{lem:general_convex_lower_bound_sqrtT} with $L=\ell$. The algorithm may be given the entire constraint function $g$. The function-value gap is bounded above using \eqref{eq:initial_gap_upper}. Condition \eqref{eq:convex_functional_accuracy} gives $\eta/r\le\epsilon/8$.

Suppose that a feasible point $\bx$ returned after at most $T$ calls is GKKT stationary. Then there exist $\lambda\ge0$, $\bar{\bg}$, and $\bar{\bw}$ such that
\[
 \br=\bar{\bg}+\lambda\bar{\bw},\qquad
 \bar{\bg}\in\partial_\delta f(\bx),\quad
 \bar{\bw}\in\partial_\delta g(\bx),\quad \|\br\|\le\epsilon.
\]
Decompose $\bar{\bg}$ as in \eqref{eq:bgj-decomp}--\eqref{eq:bgbar-decomp} and let $\beta$ be the total weight of the radial component, defined in \eqref{eq:radial_weights}. The objective sampling points lie in the full ambient ball and need not be feasible.

\par\medskip\noindent\textbf{Zone A ($\beta\le\beta^*$): Linear-Component Dominance.}\par\nobreak\smallskip
Some linear component is active at an objective sampling point. Step~4 of the proof of Lemma~\ref{lem:general_convex_lower_bound_sqrtT} gives $\|\bx\|\le2r/3+\delta<r$. Hence $g(\bx)\le-m/2$, and complementary slackness at $\bx$ itself gives
\[
 \lambda\le\frac{2\eta}{m},\qquad
 \lambda\|\bar{\bw}\|\le\frac{2\eta L_g}{m}\le\frac{\eta}{r}.
\]
Using the objective estimate of \eqref{eq:zone_a_bound},
\[
 \|\br\|\ge\frac{\ell}{12\sqrt{T+2}}-\frac{\eta}{r}
 \ge\left(\frac{4}{3}-\frac{1}{8}\right)\epsilon
 >\epsilon.
\]

\par\medskip\noindent\textbf{Zone B ($\beta>\beta^*$): Radial Activation.}\par\nobreak\smallskip
For $\bz\in\ball_\delta(\bx)$ and $\bw\in\partial g(\bz)$, convexity yields
\[
 \langle\bw,\bx\rangle
 \ge g(\bz)+m-L_g\delta.
\]
Since $\lambda|g(\bz)|\le\eta$ and $\delta\le r\le m/(2L_g)$,
\[
 \lambda\langle\bw,\bx\rangle
 \ge-\eta+\lambda(m-L_g\delta)\ge-\eta.
\]
Taking convex combinations gives
\begin{equation}\label{eq:convex_constraint_projection}
 \lambda\langle\bar{\bw},\bx\rangle\ge-\eta.
\end{equation}
If $\|\bx\|\le r$, the objective estimate \eqref{eq:zone_b_radial} gives
\[
 \langle\br,\bx\rangle
 >\frac{9\ell r}{128\sqrt{T+2}}-\eta
 \ge\left(\frac{9}{8}-\frac{1}{8}\right)r\epsilon=r\epsilon.
\]
If $\bx=\mathbf0$, this inequality contradicts $\langle\br,\mathbf0\rangle=0$. Otherwise, Cauchy--Schwarz and $\|\bx\|\le r$ give $\|\br\|>\epsilon$.

If $\|\bx\|>r$, use \eqref{eq:far_radial} and \eqref{eq:convex_constraint_projection} to obtain
\[
 \left\langle\br,\frac{\bx}{\|\bx\|}\right\rangle
 >\frac{\ell}{2}-\frac{\eta}{\|\bx\|}
 \ge\frac{\ell}{2}-\frac{\epsilon}{8}
 >\epsilon.
\]
All cases contradict $\|\br\|\le\epsilon$. This proof uses convexity, strict feasibility, and local Lipschitz bounds.
\end{proof}

\subsection{Classical precision lower bounds on a ball}\label{app:classical_precision_comparison}

\begin{proposition}[Objective-value gap estimates]\label{prop:objective_gap}
Let $f:\RR^d\to\RR$ be convex and globally $L_f$-Lipschitz.
\begin{enumerate}[(i)]
\item If $\cX=\ball_{R_{\mathrm b}}(\bc)$ and $0<\delta\le R_{\mathrm b}/2$, then
\[
f(\bx)-\min_\cX f\le3R_{\mathrm b}S_\delta(\bx)+4L_f\delta
\qquad(\bx\in\cX).
\]
\item Let $g:\RR^d\to\RR$ be convex, with $\mathcal S=\{g\le0\}\subseteq\ball_{R_{\mathrm b}}(\bc)$ and $g(\bc)=-m<0$. Assume $g$ is $L_g$-Lipschitz on a neighborhood of $\mathcal S$ containing all sampling balls, and $m>L_g\delta$. If $\bx\in\mathcal S$ is $(\delta,\epsilon,\eta)$-GKKT stationary as in Definition~\ref{def:standard_gkkt}, every associated multiplier $\lambda$ satisfies
\[
\lambda\le\frac{(L_f+\epsilon)R_{\mathrm b}+\eta}{m-L_g\delta},\qquad
f(\bx)-\min_{\mathcal S}f
\le2R_{\mathrm b}\epsilon+2L_f\delta+\eta+\lambda L_g\delta.
\]
\end{enumerate}
\end{proposition}
\begin{proof}
\textit{Part (i).} Take $\bx^*\in\operatorname*{argmin}_\cX f$ and any finite convex combination $\bv=\sum_i\omega_i(\ba_i+\bn_i)$ from Definition~\ref{def:constrained_goldstein}, with samples $\by_i\in\cX\cap\ball_\delta(\bx)$. Put $W=\sum_i\omega_i\|\bn_i\|$. The inner-ball normal inequality and $\|\ba_i\|\le L_f$ give
\[
(R_{\mathrm b}-\delta)W
\le\left\langle\sum_i\omega_i\bn_i,\bx-\bc\right\rangle
\le R_{\mathrm b}(\|\bv\|+L_f).
\]
Convexity of $f$ at each sample and $\langle\bn_i,\by_i-\bx^*\rangle\ge0$ yield
\[
\begin{aligned}
f(\bx)-f(\bx^*)
&\le2R_{\mathrm b}\|\bv\|+2L_f\delta+\delta W\\
&\le3R_{\mathrm b}\|\bv\|+4L_f\delta.
\end{aligned}
\]
Taking the infimum over $\bv$ proves (i) without requiring that the convex hull be closed.

\textit{Part (ii).} Write $\mathbf e=\bar{\ba}+\lambda\bar{\bw}$ with $\|\mathbf e\|\le\epsilon$ as in that definition. At every constraint sample $\bz$, convexity and full-ball complementary slackness give
\[
\lambda\langle\bw,\bx-\bc\rangle
\ge\lambda(g(\bz)+m-L_g\delta)
\ge\lambda(m-L_g\delta)-\eta.
\]
Average and use $\|\bar{\ba}\|\le L_f$ to bound $\lambda$. If $\bx^*$ minimizes $f$ on $\mathcal S$, the same argument at $\bx^*$ gives
\[
\lambda\langle\bar{\bw},\bx-\bx^*\rangle\ge-\eta-\lambda L_g\delta,
\qquad
f(\bx)-f(\bx^*)\le\langle\bar{\ba},\bx-\bx^*\rangle+2L_f\delta,
\]
which proves (ii).
\end{proof}

In particular, for $g(\bx)=\|\bx-\bc\|-R_{\mathrm b}$,
\[
0<\epsilon\le L_f,\quad
\delta\le\frac{R_{\mathrm b}\epsilon}{16L_f},\quad
\eta\le\frac{R_{\mathrm b}\epsilon}{8}
\quad\Longrightarrow\quad
\lambda\le\frac{34L_f}{15},\quad
f(\bx)-\min_{\mathcal S}f\le3R_{\mathrm b}\epsilon.
\]

The classical deterministic local-oracle lower bound for convex minimization on a Euclidean ball is of order $L R_{\mathrm b}/\sqrt{T+1}$ when the dimension is at least a constant multiple of $T+1$ \cite[Section 4.4.3, Proposition 1]{nemirovski1983problem}. The two estimates above therefore already give the $\epsilon^{-2}$ call-count order on balls in the stated small-neighborhood regimes. To impose the initial gap bound at $\bc$, scale the Lipschitz bound to $L=\min\{L_f,\Delta_f/R_{\mathrm b}\}$; then $f(\bc)-\min_\cX f\le L R_{\mathrm b}\le\Delta_f$. This comparison gives the precision order, rather than all the properties of the constructions in Sections~\ref{sec:framework1} and~\ref{sec:inequality_constraints}; in particular, it does not ensure an interior global minimizer.

\section{Dimension and spherical measure estimates}\label{app:exponential_estimates}

\begin{proof}[Proof of Lemma~\ref{lem:brouwer_query}]
Consider the Local Brouwer instances of Ganor, Karthik and P\'alv\"olgyi \cite[Definitions 2.1 and 5.1 and Theorem 5.2]{GKP22} on $[0,1]^n$ with the normalized Euclidean norm $\|\bw\|/\sqrt n$. At a query $\bw$, a public rule selects $k_0$ bits from each participant's private input, where $k_0$ is an absolute constant. After exchanging these $2k_0$ bits, both participants can compute $\mathbf{H}(\bw)$ and hence the next query or the output of a deterministic algorithm. The exponential communication lower bound thus gives an exponential value-query lower bound at the same constant accuracy and Lipschitz bound.

Set $C_n=[0,1/\sqrt n]^n$ and $\mathbf{P}_n=\proj_{C_n}$. For each such $\mathbf{H}$, define
\[
\mathbf{F}_n(\bz)=\mathbf{H}(\sqrt n \mathbf{P}_n(\bz))/\sqrt n
\]
on $\ball_1^n(\mathbf{0})$. The map $\mathbf{F}_n$ takes values in $C_n$ and has the same Lipschitz bound as $\mathbf{H}$ by nonexpansiveness of $\mathbf{P}_n$. Each value of $\mathbf{F}_n$ uses one value of $\mathbf{H}$. Moreover, at $\bw=\sqrt n \mathbf{P}_n(\bz)$,
\[
\frac{\|\bw-\mathbf{H}(\bw)\|}{\sqrt n}
=\|\mathbf{P}_n(\bz)-\mathbf{P}_n(\mathbf{F}_n(\bz))\|
\le\|\bz-\mathbf{F}_n(\bz)\|.
\]
This transfers the lower bound from the scaled cube to the unit ball in the source dimension.

The source construction has $n=\Theta(\log N)$ for private inputs of length $N$. For a sufficiently large $d$, take $N=2^{\lfloor d/c_1\rfloor}$ with an absolute positive integer $c_1$ large enough to ensure $n\le d$ and $n=\Omega(d)$. Extend $\mathbf{F}_n$ to $\ball_1^d(\mathbf{0})$ by $\bz\mapsto(\mathbf{F}_n(\bz_{1:n}),\mathbf{0}_{d-n})$, where $\bz_{1:n}$ consists of the first $n$ coordinates. The extension has the same Lipschitz bound and uses one value of $\mathbf{F}_n$ per query. Taking the first $n$ coordinates of an approximate fixed point of the extension gives an approximate fixed point of $\mathbf{F}_n$ at the same accuracy. This proves \eqref{eq:brouwer_query} after adjusting the absolute constants.
\end{proof}

\paragraph{Spherical caps and product cones.}
For $0<b\le1$ and $\bv\in\mathbb S^{m-1}$, the chordal cap
\[
A:=\{\bu\in\mathbb S^{m-1}\mid\|\bu-\bv\|\le b\}
\]
satisfies $\mathfrak m_m(A)\ge(b/8)^m$. For $m\ge2$, project the spherical graph onto the tangent hyperplane at its pole. The graph over the disk of radius $b/2$ lies in the cap. Let $V_m$ be the unit ball volume, so the unit sphere has area $mV_m$. Slicing the ball gives $V_m\le2V_{m-1}$, hence
\[
\mathfrak m_m(A)
\ge\frac{V_{m-1}(b/2)^{m-1}}{mV_m}
\ge\frac{(b/2)^{m-1}}{2m}\ge\left(\frac b8\right)^m.
\]
For $m=1$, the sphere is $\{-1,1\}$ with probability $1/2$ at each point. Independent standard Gaussian blocks have independent uniform directions. For the product cone in Lemma~\ref{lem:constraint_cones}, using $b=1/(4K)$ in the constrained blocks and probability $1$ in the remaining blocks yields
\[
\mathfrak m_d\left(\left(\prod_{j=1}^k C_j\right)\cap\unit\right)
\ge\prod_{j=1}^k\left(\frac1{32K}\right)^{d_j}=\kappa^d.
\]

To bound the measure of a pair of caps at opposite poles, use the same graph representation. If the projected disks have radius $\alpha\le1/4$, the graph area element is bounded above by twice the area element on the disk. The sphere has area at least $2V_{d-1}$, since its two hemispheres each project onto a unit ball in dimension $d-1$. Thus the total normalized measure of the two caps is at most $2\alpha^{d-1}$.

\paragraph{Directions forming an acute angle with $\bp$.}
Let $C$ be the cone in Lemma~\ref{lem:constraint_cones} and let $\bp\in C\cap\unit$.
For the spherical directions making an acute angle with $\bp$, the required estimate is
\[
\mathfrak m_d(\{\bu\in C\cap\unit\mid\langle\bu,\bp\rangle\ge1/\sqrt2\})
\ge(\kappa/2)^d.
\]
 Reflection of $\bv\in C$ with $\langle \bv,\bp\rangle<0$ across $\bp^\perp$ adds a nonnegative multiple of $\bp$; therefore at least half the cone's measure lies in the hemisphere $\langle \bv,\bp\rangle\ge0$. The map $\bv\mapsto(\bp+\bv)/\|\bp+\bv\|$ stays in the cone and halves the polar angle in that hemisphere. Its spherical Jacobian is $\frac12(\sin(t/2)/\sin t)^{d-2}\ge2^{-(d-1)}$ for polar angles $0<t\le\pi/2$.

\begin{proof}[Proof of Corollary~\ref{cor:exp_fixed_accuracy}]
Choose
\[
 \epsilon_0=\min\left\{\frac{\beta_{\mathrm B}}{8},\frac{\kappa}{128}\right\},
 \qquad c_0=\min\left\{\frac{b_{\mathrm B}}{4},\frac{\log 2}{8}\right\}.
\]
For all sufficiently large $d$, $\exp(c_0d)\ge2$ and
\[
 T+2\le2\exp(c_0d)\le\exp(b_{\mathrm B}d).
\]
Since $T+2$ is an integer, $T+2\le Q_{\mathrm B}(d)$. The inequality $32\epsilon/L_f\le\kappa/4$ implies that the ratio of the two sides of the angular inequality in \eqref{eq:exp_parameters} is at most
\[
 \frac{16}{\kappa}\exp\bigl(-(\log 2-2c_0)d\bigr).
\]
This is less than one for all sufficiently large $d$, which gives the result.
\end{proof}

\section{Comparison of constrained stationarity measures}\label{app:ggsp_comparison}

\begin{proof}[Proof of Proposition~\ref{prop:ggsp_comparison}]
\textit{Parts (i) and (ii).} Suppose $N_{\cX}(\by)\subseteq N_{\cX}(\bx)$ for every $\by\in \ball_\delta(\bx)\cap \cX$. Each vector in the convex hull in Definition~\ref{def:constrained_goldstein} can then be written as $\bv+\bn$, with $\bv\in\partial_\delta f(\bx)$ and $\bn\in N_{\cX}(\bx)$. Since $\proj_{\cX}(\bx+\gamma \bn)=\bx$, nonexpansiveness gives $\|\bw_{\bv}\|\le\|\bv+\bn\|$, and taking the infimum proves $M_{\gamma,\delta}(\bx)\le S_\delta(\bx)$.

For part (ii), suppose $\ball_{r_{\mathrm{in}}}(\bc)\subseteq \cX$, let $0<\delta<r_{\mathrm{in}}$, and assume $f$ is globally $L_f$-Lipschitz. We show that
\[
S_\delta(\bx)\le\epsilon
\quad\Longrightarrow\quad
M_{\gamma,\delta}(\bx)\le\epsilon+
\sqrt{\frac{\delta\,\|\bx-\bc\|(L_f+\epsilon)}
{\gamma(r_{\mathrm{in}}-\delta)}}.
\]
Take $\epsilon'>\epsilon$ and choose a convex combination $\mathbf{e}=\sum_i a_i(\bv_i+\bn_i)$ of feasible samples with $\|\mathbf{e}\|\le\epsilon'$. Here $\bv_i\in\partial f(\by_i)$, $\bn_i\in N_{\cX}(\by_i)$, and $\by_i\in \cX\cap \ball_\delta(\bx)$. Let $\bv=\sum_i a_i\bv_i$ and $\bn=\sum_i a_i\bn_i$. Then $\|\bv\|\le L_f$ and $\|\bn\|\le L_f+\epsilon'$. The inner ball gives $\langle \bn_i,\by_i-\bc\rangle\ge r_{\mathrm{in}}\|\bn_i\|$, hence
\[
(r_{\mathrm{in}}-\delta)\sum_i a_i\|\bn_i\|
\le\langle \bn,\bx-\bc\rangle\le\|\bx-\bc\|(L_f+\epsilon').
\]
For $\bp_{\bv}$ in \eqref{eq:ggsp_projection}, feasibility gives $\langle \bn_i,\bp_{\bv}-\bx\rangle\le\langle \bn_i,\by_i-\bx\rangle\le\delta\|\bn_i\|$. Also $\langle \bn_{\bv},\bw_{\bv}\rangle\le0$. Combining these inequalities gives
\[
\begin{aligned}
\|\bw_{\bv}\|^2
&\le\langle \bw_{\bv},\bv\rangle
=\langle \bw_{\bv},\mathbf{e}\rangle+\frac{\langle \bn,\bp_{\bv}-\bx\rangle}{\gamma}\\
&\le\epsilon'\|\bw_{\bv}\|
+\frac{\delta\,\|\bx-\bc\|(L_f+\epsilon')}{\gamma(r_{\mathrm{in}}-\delta)}.
\end{aligned}
\]
Solving the quadratic inequality and letting $\epsilon'\downarrow\epsilon$ gives part (ii).

\textit{Part (iii).} Suppose $\bv$ is a convex combination of subgradients at feasible samples in $\ball_\delta(\bx)$ and satisfies $\|\bw_{\bv}\|\le\epsilon$. Let $R_\gamma=\max\{\delta,\gamma\epsilon\}$. Then $\bp_{\bv}\in \cX\cap \ball_{R_\gamma}(\bx)$ and
\[
S_{R_\gamma}(\bx)\le\epsilon.
\]
Indeed, for $\bv_{\bp}\in\partial f(\bp_{\bv})$ and $t\in(0,1)$, the vector $(1-t)\bv+t(\bv_{\bp}+\bn_{\bv}/t)$ lies in the convex hull defining $S_{R_\gamma}(\bx)$, using zero normals for the samples that form $\bv$. As $t\downarrow0$ this vector converges to $\bv+\bn_{\bv}=\bw_{\bv}$, whose norm is at most $\epsilon$. If $\ball_\delta(\bx)\subseteq \cX$, every objective sampling point is feasible. If also $\gamma\epsilon\le\delta$, then the same argument gives $S_\delta(\bx)\le\epsilon$ with the original radius.

\textit{Part (iv).} Suppose $f$ is continuously differentiable on a neighborhood of $\ball_{R_\gamma}(\bx)$ and has $L_{\mathrm s}$-Lipschitz gradient on that ball. Choose $\bv\in\partial_\delta f(\bx)$ with $\|\bw_{\bv}\|\le\epsilon$. The gradient Lipschitz bound gives $\|\bv-\nabla f(\bp_{\bv})\|\le L_{\mathrm s}(\delta+\gamma\epsilon)$. Applying \eqref{eq:ggsp_projection} at $\bp_{\bv}$ and using $\bp_{\bv}\in \ball_{\gamma\epsilon}(\bx)\cap \cX$ yields
\[
S_{\gamma\epsilon}(\bx)\le\epsilon+L_{\mathrm s}(\delta+\gamma\epsilon).
\]

\end{proof}

\begin{proof}[Proof of Proposition~\ref{prop:residual_comparison}]
For part (i), the assumed normal-cone inclusion and convexity give
\[
\operatorname{conv}\bigcup_{\by\in\ball_\delta(\bx)\cap\cX}
       \bigl(\partial f(\by)+N_\cX(\by)\bigr)
\subseteq\partial_\delta f(\bx)+N_\cX(\bx).
\]
Taking distances from the origin proves the inequality.

For part (ii), the ball inclusion implies $\bx\in\operatorname{int}(\cX)$, and hence $N_\cX(\bx)=\{\mathbf0\}$ and $\mathcal R_\delta(\bx)=\dist(\mathbf0,\partial_\delta f(\bx))$. Every objective sampling point is feasible. Choosing zero normals at these points gives $\partial_\delta f(\bx)\subseteq\partial_\delta^\cX f(\bx)$ and hence $S_\delta(\bx)\le\mathcal R_\delta(\bx)$. If the entire ball lies in $\operatorname{int}(\cX)$, all sampled normals vanish and the two sets coincide.

For part (iii), choose $\bv\in\partial_\delta f(\bx)$ and $\bn\in N_\cX(\bx)$ with $\|\bv+\bn\|=\mathcal R_\delta(\bx)$. Such a pair exists by compactness of $\partial_\delta f(\bx)$ and closedness of the normal cone. The gradient Lipschitz bound gives $\|\bv-\nabla f(\bx)\|\le L_{\mathrm s}\delta$. Since $\bx$ itself is an admissible sample,
\[
S_\delta(\bx)
\le\|\nabla f(\bx)+\bn\|
\le\|\bv+\bn\|+L_{\mathrm s}\delta
=\mathcal R_\delta(\bx)+L_{\mathrm s}\delta.
\]
\end{proof}

\begin{acknowledgements}
The authors thank Luo Luo for suggesting a comparison with the stationarity criterion in \cite{zero}, and Guy Kornowski for helpful comments on the relationship to existing lower bounds and for pointing out the related work \cite{KFO26,BK26}.
\end{acknowledgements}

\section*{Statement on AI use}
The authors used ChatGPT and Codex (OpenAI), as well as Gemini (Google), to assist with mathematical arguments, literature review, and manuscript preparation. The authors retain full responsibility for the mathematical content, references, and presentation of the manuscript.

\section*{Statements and Declarations}
\paragraph{Funding.}
This work was partially supported by the National Natural Science Foundation of China (No. 12288201).

\paragraph{Competing interests.}
The authors have no competing interests to declare that are relevant to the content of this article.

\paragraph{Data availability.}
No data are associated with this article. Data sharing is not applicable to this article.

\paragraph{Author contributions.}
Yixi Ding: Formal analysis, Investigation, Writing -- original draft. Ya-xiang Yuan: Supervision.

\bibliographystyle{spmpsci}
\bibliography{references}
\addcontentsline{toc}{section}{References}
\end{document}